\documentclass[10pt, a4paper]{amsart}
\usepackage{amsfonts}
\usepackage{amssymb, amsmath}
\usepackage{bbm}
\usepackage{url}

\newtheorem{theorem}{Theorem}[section]
\newtheorem{proposition}[theorem]{Proposition}
\newtheorem{corollary}[theorem]{Corollary}

\theoremstyle{definition}
\newtheorem{definition}[theorem]{Definition}
\newtheorem{example}[theorem]{\textbf{Example}}
\newtheorem{remark}[theorem]{\textbf{Remark}}

\newtheorem{notation}[theorem]{\textbf{Notation}}

\title[Cohen-Macaulayness of tangent cones of monomial curves in $\mathbb{A}^{4}(K)$]{On the Cohen-Macaulayness of tangent cones of monomial curves in $\mathbb{A}^{4}(K)$}

\author[F. Arslan]{Feza Arslan}
\address {Department of Mathematics, Mimar Sinan Fine Arts University, Istanbul, 34427, Turkey} \email{sefa.feza.arslan@msgsu.edu.tr}
\author[A. Katsabekis]{Anargyros Katsabekis}
\address {Department of Mathematics, Mimar Sinan Fine Arts University, Istanbul, 34427, Turkey} \email{katsabek@aegean.gr}
\author[M. Nalbandiyan]{ Melissa Nalbandiyan}
\address {Department of Mathematics, Mimar Sinan Fine Arts University, Istanbul, 34427, Turkey} \email{assilem89@gmail.com}

\thanks{The second author has been supported by TUBITAK 2221 Visiting Scientists and Scientists on Sabbatical Leave Fellowship Program}

\begin{document}

\date{\today}

\begin{abstract} In this paper we give necessary and sufficient conditions for the Cohen-Macaulayness of the tangent cone of a monomial curve in the 4-dimensional affine space. We study particularly the case where $C$ is a Gorenstein non-complete intersection monomial curve.

\end{abstract}

\maketitle

\section{Introduction}

Cohen-Macaulayness of tangent cones of monomial curves has been studied by many authors, see for instance \cite{Ar}, \cite{ArMe}, \cite{Ga}, \cite{Her1}, \cite{Her2}, \cite{Ja}, \cite{Mo2}, \cite{Mo1}, \cite{Pa}, \cite{Ro}, \cite{Sa}, \cite{Sh}, \cite{She}. It constitutes an important problem, since for example Cohen-Macaulayness of the tangent cone guarantees that the Hilbert function of the local ring associated to the monomial curve is non-decreasing and therefore reduces its computation to the computation of the Hilbert function of an Artin local ring.

In this article, we first deal with the above problem in the case of any monomial curve in the 4-dimensional affine space $\mathbb{A}^{4}(K)$, where $K$ is a field. In section 2, by using the classification in terms of critical binomials given by Katsabekis and Ojeda \cite{KaOj}, we study in detail the problem for Case 1 in this classification and give sufficient conditions for the Cohen-Macaulayness of the tangent cone. We consider the remaining cases in the Appendix, where we give in all these cases necessary and sufficient conditions for the Cohen-Macaulayness of the tangent cone. 

In section 3 we consider the problem for non-complete intersection Gorenstein monomial curves. In this case, Bresinsky has not only shown that there is a minimal generating set for the defining ideal of the monomial curve consisting of five generators, but also given the explicit form of these generators \cite{Bresinsky75}. Actually there are 6 permutations of the above generator set. It is worth to note that Theorem 2.10 in \cite{ArMe} provides a sufficient condition for the Cohen-Macaulayness of the tangent cone in four of the aforementioned cases. In this paper, we generalize their result and provide a necessary and sufficient condition for the Cohen-Macaulayness of the tangent cone in all 6 six permutations. Finally, we use these results to give some families of Gorenstein monomial curves in $\mathbb{A}^{4}(K)$ with corresponding local rings having non-decreasing Hilbert function, thus giving a partial answer to Rossi's problem \cite{Rossi}. This problem asks whether the Hilbert function of a Gorenstein local ring of dimension one is non-decreasing. Recently, it has been shown that there are many families of monomial curves giving negative answer to this problem \cite{OST}, but one should note that Rossi's problem is still open for Gorenstein local rings associated to monomial curves in $\mathbb{A}^{4}(K)$. 

Let $\{n_{1},\ldots,n_{d}\}$ be a set of all-different positive
integers with ${\rm gcd}(n_{1},\ldots,n_{d})=1$. Let $K[x_{1},\ldots,x_{d}]$ be the polynomial ring in $d$ variables. We shall denote by ${\bf x}^{\bf u}$ the monomial $x_{1}^{u_1} \cdots x_{d}^{u_d}$ of $K[x_{1},\ldots,x_{d}]$, with ${\bf u}=(u_{1},\ldots,u_{d}) \in \mathbb{N}^{d}$, where $\mathbb{N}$ stands for the set of non-negative integers. Consider the affine monomial curve in the $d$-space $\mathbb{A}^{d}(K)$ defined parametrically by $$x_{1}=t^{n_1},\ldots,x_{d}=t^{n_d}.$$ The toric ideal of $C$, denoted by $I(C)$, is the kernel of the $K$-algebra homomorphism $\phi:K[x_{1},\ldots,x_{d}] \rightarrow K[t]$ given by $$\phi(x_{i})=t^{n_i} \ \ \textrm{for all} \ \ 1 \leq i \leq d.$$

We grade $K[x_{1},\ldots,x_{d}]$ by the semigroup $\mathcal{S}:=\{g_{1}n_{1}+\cdots+g_{d}n_{d}|g_{i} \in \mathbb{N}\}$ setting ${\rm deg}_{\mathcal{S}}(x_{i})=n_{i}$ for $i=1,\ldots,d$. The $\mathcal{S}$-degree of a monomial ${\bf x}^{\bf u}=x_{1}^{u_{1}} \cdots x_{d}^{u_{d}}$ is defined by $${\rm deg}_{\mathcal{S}}({\bf x}^{\bf u})=u_{1}n_{1}+\cdots+u_{d}n_{d} \in \mathcal{S}.$$ The ideal $I(C)$ is generated by all the binomials ${\bf x}^{{\bf u}}-{\bf x}^{{\bf v}}$ such that ${\rm deg}_{\mathcal{S}}({\bf x}^{{\bf u}})={\rm deg}_{\mathcal{S}}({\bf x}^{{\bf v}})$ see for example, \cite[Lemma 4.1]{Sturmfels95}.

Let $n_{1}<n_{2}<\cdots<n_{d}$ and $n_{1}+\mathcal{S}=\{n_{1}+m|m \in \mathcal{S}\}$.

\begin{theorem} \label{BasicHer} (\cite{Her1}) The monomial curve $C$ has Cohen-Macaulay tangent cone at the origin if and only if for all integers $v_{2} \geq 0, v_{3} \geq 0,\ldots,v_{d} \geq 0$ such that $\sum_{i=2}^{d}v_{i}n_{i} \in n_{1}+\mathcal{S}$, there exist $w_{1}>0$, $w_{2} \geq 0, \ldots,w_{d} \geq 0$ such that $\sum_{i=2}^{d}v_{i}n_{i}=\sum_{i=1}^{d}w_{i}n_{i}$ and $\sum_{i=2}^{d}v_{i} \leq \sum_{i=1}^{d}w_{i}$.

\end{theorem}

Note that $x_{1}^{\frac{n_i}{{\rm gcd}(n_{1},n_{i})}}-x_{i}^{\frac{n_1}{{\rm gcd}(n_{1},n_{i})}} \in I(C)$ and also $\frac{n_i}{{\rm gcd}(n_{1},n_{i})}>\frac{n_1}{{\rm gcd}(n_{1},n_{i})}$, for every $2 \leq i \leq d$. Thus, to decide the Cohen-Macaulayness of the tangent cone of $C$ it suffices to consider only such $v_i$ with the extra condition that $v_{i}<\frac{n_1}{{\rm gcd}(n_{1},n_{i})}$.

The computations of this paper are performed by using CoCoA (\cite{Coc}).

\section{The general case}

Let $A=\{n_{1},\ldots,n_{4}\}$ be a set of relatively prime positive integers.

\begin{definition}
{\rm A binomial $x_i^{a_i} - \prod_{j \neq i} x_j^{u_{ij}} \in I(C)$ is called \textbf{critical} with respect to $x_i$ if $a_i$ is the least positive integer such that $a_i n_i \in \sum_{j \neq i} \mathbb{N} n_j.$
The \textbf{critical ideal} of $A$, denoted by $\mathcal{C}_{A}$, is the ideal of $K[x_{1},\ldots,x_{4}]$ generated by all the critical binomials of $I(C).$}
\end{definition}

The support ${\rm supp}({\bf x}^{\bf u})$ of a monomial ${\bf x}^{\bf u}$ is the set $${\rm supp}({\bf x}^{\bf u})=\{i \in \{1,\ldots,4\}|x_{i} \ \textrm{divides} \ {\bf x}^{\bf u}\}.$$ The support of a binomial $B={\bf x}^{\bf u}-{\bf x}^{\bf v}$ is the set ${\rm supp}({\bf x}^{\bf u}) \cup {\rm supp}({\bf x}^{\bf v})$. If the support of $B$ equals the set $\{1,\ldots,4\}$, then we say that $B$ has full support. Let $\mu(\mathcal{C}_{A})$ be the minimal number of generators of the ideal $\mathcal{C}_{A}$.

\begin{theorem} (\cite{KaOj}) \label{KaOjVeryBasic} After permuting the variables, if necessary, there exists a minimal
system of binomial generators $S$ of the critical ideal $\mathcal{C}_{A}$ of the following form: \begin{itemize}
\item[CASE 1:] If $a_i n_i \neq a_j n_j$, for every $i \neq j,$ then $S = \{ x_i^{a_i} -
{\bf x}^{{\bf u}_i},\ i = 1, \ldots, 4 \}$ \item[CASE 2:] If $a_1 n_1= a_2 n_2$ and $a_3 n_3 = a_4 n_4,$
then either $a_2 n_2 \neq a_3 n_3$ and \begin{itemize}  \item[(a)] $S= \{x_1^{a_1} - x_2^{a_2}, x_3^{a_3} - x_4^{a_4},
x_4^{a_4} -{\bf x}^{{\bf u}_4} \}$ when $\mu(\mathcal{C}_{A})=3$ \item[(b)] $S= \{x_1^{a_1} -
x_2^{a_2}, x_3^{a_3} - x_4^{a_4}\}$ when $\mu(\mathcal{C}_{A})=2$ \end{itemize} or $a_2 n_2 =a_3 n_3$ and
\begin{itemize} \item[(c)] $S = \{x_1^{a_1} - x_2^{a_2}, x_2^{a_2} - x_3^{a_3}, x_3^{a_3} -
x_4^{a_4}\}$ \end{itemize} \smallskip \item[CASE 3:] If $a_1 n_1 = a_2 n_2 = a_3 n_3 \neq a_4 n_4,$ then
$S= \{ x_1^{a_1} - x_2^{a_2}, x_2^{a_2} - x_3^{a_3}, x_4^{a_4} - {\bf x}^{{\bf u}_4}\}$
\item[CASE 4:] If $a_1 n_1 = a_2 n_2$ and $a_i n_i \neq a_j n_j$ for all $\{i, j\} \neq \{1, 2\},$ then
\begin{itemize} \item[(a)] $S= \{x_1^{a_1} - x_2^{a_2}, x_i^{a_i} - {\bf x}^{{\bf u}_i}\ \mid\ i =
2,3,4\}$ when $\mu(\mathcal{C}_{A})=4$ \item[(b)] $S= \{x_1^{a_1} - x_2^{a_2}, x_i^{a_i} -
{\bf x}^{{\bf u}_i}\ \mid\ i = 3,4\}$ when $\mu(\mathcal{C}_{A})=3$ \end{itemize} \end{itemize} where, in each case, ${\bf x}^{{\bf u}_i}$ denotes an appropriate monomial whose support has cardinality greater than or equal to two.
\end{theorem}

\begin{theorem} (\cite{KaOj}) \label{Katsabekisbasic}
The union of $S,$ the set $\mathcal{I}$ of all binomials $x_{i_1}^{u_{i_1}} x_{i_2}^{u_{i_2}} - x_{i_3}^{u_{i_3}} x_{i_4}^{u_{i_4}} \in I(C)$ with $0 < u_{i_j} < a_j,\ j = 1,2$, $u_{i_3}>0$, $u_{i_4}>0$ and $x_{i_3}^{u_{i_3}} x_{i_4}^{u_{i_4}}$ indispensable, and the set $\mathcal{R}$ of all binomials $x_1^{u_1} x_2^{u_2} - x_3^{u_3} x_4^{u_4} \in I(C) \setminus \mathcal{I}$ with full support such that
\begin{itemize}
\item $u_1 \leq a_1$ and $x_3^{u_3} x_4^{u_4}$ is indispensable, in the CASES 2(a) and 4(b).
\item $u_1 \leq a_1$ and/or $u_3 \leq a_3$ and there is no $x_1^{v_1} x_2^{v_2} - x_3^{v_3} x_4^{v_4} \in I(C)$ with full support such that $x_1^{v_1} x_2^{v_2}$ properly divides $x_1^{u_1 + \alpha c_1} x_2^{u_2 - \alpha c_2}$ or $x_3^{v_3} x_4^{v_4}$ properly divides $x_3^{u_3 + \alpha c_3} x_4^{u_4 - \alpha u_4}$ for some $\alpha \in \mathbb{N},$ in the CASE 2(b).
\end{itemize}
is a minimal system of generators of $I(C)$ (up to permutation of indices).
\end{theorem}

A binomial $B \in I(C)$ is called indispensable of $I(C)$ if every system of binomial generators of $I(C)$ contains $B$ or $-B$. By Corollary 2.16 in \cite{KaOj} every $f \in \mathcal{I}$ is an indispensable binomial of $I(C)$.

\begin{notation} {\rm Given a monomial ${\bf x}^{\bf u}$ we will write ${\rm deg}({\bf x}^{\bf u}):=\sum_{i=1}^{4}u_{i}$.}
\end{notation}

For the rest of this section we will assume that $n_{1}<n_{2}<n_{3}<n_{4}$. To prove our results we will make repeated use of Theorem \ref{BasicHer}.

\begin{proposition} \label{1(b)} Suppose that $I(C)$ is given as in case 1. Let $S=\{x_1^{a_1} -
{\bf x}^{{\bf u}}, x_2^{a_2} -
{\bf x}^{{\bf v}}, x_3^{a_3} -
{\bf x}^{{\bf w}}, x_4^{a_4} -
{\bf x}^{{\bf z}}\}$ be a generating set of $\mathcal{C}_{A}$ and let $1 \in {\rm supp}({\bf x}^{{\bf v}})$. In the following cases $C$ has Cohen-Macaulay tangent cone at the origin. (i) $1 \in {\rm supp}({\bf x}^{{\bf w}})$, $1 \in {\rm supp}({\bf x}^{{\bf z}})$ and  \begin{enumerate} \item $a_{2} \leq {\rm deg}({\bf x}^{{\bf v}})$,\item $a_{3} \leq {\rm deg}({\bf x}^{{\bf w}})$ and \item for every binomial $f=M-N \in \mathcal{I}$ with $1 \in {\rm supp}(M)$ we have that ${\rm deg}(N) \leq {\rm deg}(M)$.
\end{enumerate}
(ii) $1 \in {\rm supp}({\bf x}^{{\bf w}})$, $1 \notin {\rm supp}({\bf x}^{{\bf z}})$ and \begin{enumerate} \item $a_{2} \leq {\rm deg}({\bf x}^{{\bf v}})$, \item $a_{3} \leq {\rm deg}({\bf x}^{{\bf w}})$ and \item for every monomial $M=x_{2}^{d_{2}}x_{3}^{d_{3}}x_{4}^{d_4}$, where $d_{2}<a_{2}$ and $d_{3}<a_{3}$, with $d_{2}n_{2}+d_{3}n_{3}+d_{4}n_{4} \in n_{1}+\mathcal{S}$ there exists a monomial $N$ with $1 \in {\rm supp}(N)$ such that $M-N \in I(C)$ and also ${\rm deg}(M) \leq {\rm deg}(N)$.
\end{enumerate}
(iii) $1 \notin {\rm supp}({\bf x}^{{\bf w}})$, $1 \in {\rm supp}({\bf x}^{{\bf z}})$ and \begin{enumerate} \item $a_{2} \leq {\rm deg}({\bf x}^{{\bf v}})$, \item for every monomial $M=x_{2}^{d_{2}}x_{3}^{d_{3}}x_{4}^{d_{4}}$, where $d_{2}<a_{2}$ and $d_{4}<a_{4}$, with $d_{2}n_{2}+d_{3}n_{3}+d_{4}n_{4} \in n_{1}+\mathcal{S}$ there exists a monomial $N$ with $1 \in {\rm supp}(N)$ such that $M-N \in I(C)$ and also ${\rm deg}(M) \leq {\rm deg}(N)$.
\end{enumerate}
(iv) $1 \notin {\rm supp}({\bf x}^{{\bf w}})$, $1 \notin {\rm supp}({\bf x}^{{\bf z}})$ and \begin{enumerate} \item $a_{2} \leq {\rm deg}({\bf x}^{{\bf v}})$, \item for every monomial $M=x_{2}^{d_{2}}x_{3}^{d_{3}}x_{4}^{d_{4}}$, where $d_{2}<a_{2}$, with $d_{2}n_{2}+d_{3}n_{3}+d_{4}n_{4} \in n_{1}+\mathcal{S}$ there exists a monomial $N$ with $1 \in {\rm supp}(N)$ such that $M-N \in I(C)$ and also ${\rm deg}(M) \leq {\rm deg}(N)$.
\end{enumerate}
\end{proposition}

\noindent \textbf{Proof.} (i) Let $M=x_{2}^{d_2}x_{3}^{d_3}x_{4}^{d_4}$, where $d_{2} \geq 0$, $d_{3} \geq 0$ and $d_{4} \geq 0$, with $d_{2}n_{2}+d_{3}n_{3}+d_{4}n_{4} \in n_{1}+\mathcal{S}$. Thus there exists a monomial $P$ such that $1 \in {\rm supp}(P)$ and $M-P \in I(C)$. Let $d_{2} \geq a_{2}$, then we consider the monomial $P=x_{2}^{d_{2}-a_{2}}x_{3}^{d_3}x_{4}^{d_4}{\bf x}^{{\bf v}}$. We have that $M-P \in I(C)$ and also $${\rm deg}(M)=d_{2}+d_{3}+d_{4} \leq d_{2}+d_{3}+d_{4}+({\rm deg}({\bf x}^{{\bf v}})-a_{2})={\rm deg}(P).$$ Let $d_{3} \geq a_{3}$, then we consider the monomial $P=x_{2}^{d_{2}}x_{3}^{d_{3}-a_{3}}x_{4}^{d_4}{\bf x}^{{\bf w}}$. We have that $M-P \in I(C)$ and also $${\rm deg}(M)=d_{2}+d_{3}+d_{4} \leq d_{2}+d_{3}+d_{4}+({\rm deg}({\bf x}^{{\bf w}})-a_{3})={\rm deg}(P).$$ Let $d_{4} \geq a_{4}$, then we consider the monomial $P=x_{2}^{d_{2}}x_{3}^{d_{3}}x_{4}^{d_{4}-a_{4}}{\bf x}^{{\bf z}}$. We have that $M-P \in I(C)$ and also $${\rm deg}(M)=d_{2}+d_{3}+d_{4}<d_{2}+d_{3}+d_{4}+({\rm deg}({\bf x}^{{\bf z}})-a_{4})={\rm deg}(P).$$ Suppose that $d_{2}<a_{2}$, $d_{3}<a_{3}$ and $d_{4}<a_{4}$. There are 2 cases: (1) there exists a binomial $f=G-H \in \mathcal{I}$ with $1 \in {\rm supp}(G)$ such that $H$ divides $M$, so $M=HH'$. Note that ${\rm deg}(H) \leq {\rm deg}(G)$ from condition (3). Then $GH'-HH' \in I(C)$ and also $${\rm deg}(M)={\rm deg}(HH') \leq {\rm deg}(G)+{\rm deg}(H').$$
(2) there exists no binomial $f=G-H$ in $\mathcal{I}$ such that $H$ divides $M$. Recall that $M-P \in I(C)$, $S \cup \mathcal{I}$ generates $I(C)$ and also $1 \in {\rm supp}({\bf x}^{{\bf v}})$, $1 \in {\rm supp}({\bf x}^{{\bf w}})$ and $1 \in {\rm supp}({\bf x}^{{\bf z}})$. Then necessarily ${\bf x}^{{\bf u}}$ divides $M$, so $M={\bf x}^{{\bf u}}M'$. Let $P=M'x_{1}^{a_{1}}$, then $M-P \in I(C)$ and also $${\rm deg}(M)={\rm deg}({\bf x}^{{\bf u}})+{\rm deg}(M')<{\rm deg}(x_{1}^{a_1})+{\rm deg}(M')={\rm deg}(P).$$
(ii) Let $M=x_{2}^{d_2}x_{3}^{d_3}x_{4}^{d_4}$, where $d_{2} \geq 0$, $d_{3} \geq 0$ and $d_{4} \geq 0$, with $d_{2}n_{2}+d_{3}n_{3}+d_{4}n_{4} \in n_{1}+\mathcal{S}$. Let $d_{2} \geq a_{2}$, then we consider the monomial $P=x_{2}^{d_{2}-a_{2}}x_{3}^{d_3}x_{4}^{d_4}{\bf x}^{{\bf v}}$. We have that $M-P \in I(C)$ and also $${\rm deg}(M)=d_{2}+d_{3}+d_{4} \leq d_{2}+d_{3}+d_{4}+({\rm deg}({\bf x}^{{\bf v}})-a_{2})={\rm deg}(P).$$ Let $d_{3} \geq a_{3}$, then we consider the monomial $P=x_{2}^{d_{2}}x_{3}^{d_{3}-a_{3}}x_{4}^{d_4}{\bf x}^{{\bf w}}$. We have that $M-P \in I(C)$ and also $${\rm deg}(M)=d_{2}+d_{3}+d_{4} \leq d_{2}+d_{3}+d_{4}+({\rm deg}({\bf x}^{{\bf w}})-a_{3})={\rm deg}(P).$$ Suppose that $d_{2}<a_{2}$ and $d_{3}<a_{3}$, then from condition (3) we are done.\\
(iii) Let $M=x_{2}^{d_2}x_{3}^{d_3}x_{4}^{d_4}$, where $d_{2} \geq 0$, $d_{3} \geq 0$ and $d_{4} \geq 0$, with $d_{2}n_{2}+d_{3}n_{3}+d_{4}n_{4} \in n_{1}+\mathcal{S}$. Let $d_{2} \geq a_{2}$, then we consider the monomial $P=x_{2}^{d_{2}-a_{2}}x_{3}^{d_3}x_{4}^{d_4}{\bf x}^{{\bf v}}$. We have that $M-P \in I(C)$ and also $${\rm deg}(M)=d_{2}+d_{3}+d_{4} \leq d_{2}+d_{3}+d_{4}+({\rm deg}({\bf x}^{{\bf v}})-a_{2})={\rm deg}(P).$$ Let $d_{4} \geq a_{4}$, then we consider the monomial $P=x_{2}^{d_{2}}x_{3}^{d_{3}}x_{4}^{d_{4}-a_{4}}{\bf x}^{{\bf z}}$. We have that $M-P \in I(C)$ and also $${\rm deg}(M)=d_{2}+d_{3}+d_{4}<d_{2}+d_{3}+d_{4}+{\rm deg}({\bf x}^{{\bf z}})-a_{4}={\rm deg}(P).$$ Suppose that $d_{2}<a_{2}$ and $d_{4}<a_{4}$, then from condition (2) we are done.\\
(iv) Let $M=x_{2}^{d_2}x_{3}^{d_3}x_{4}^{d_4}$, where $d_{2} \geq 0$, $d_{3} \geq 0$ and $d_{4} \geq 0$, with $d_{2}n_{2}+d_{3}n_{3}+d_{4}n_{4} \in n_{1}+\mathcal{S}$. Let $d_{2} \geq a_{2}$, then we consider the monomial $P=x_{2}^{d_{2}-a_{2}}x_{3}^{d_3}x_{4}^{d_4}{\bf x}^{{\bf v}}$. We have that $M-P \in I(C)$ and also $${\rm deg}(M)=d_{2}+d_{3}+d_{4} \leq d_{2}+d_{3}+d_{4}+({\rm deg}({\bf x}^{{\bf v}})-a_{2})={\rm deg}(P).$$ Suppose that $d_{2}<a_{2}$, then from condition (2) we are done.

\begin{proposition} \label{2(a)} Suppose that $I(C)$ is given as in case 1. Let $S=\{x_1^{a_1} -
{\bf x}^{{\bf u}}, x_2^{a_2} -
{\bf x}^{{\bf v}}, x_3^{a_3} -
{\bf x}^{{\bf w}}, x_4^{a_4} -
{\bf x}^{{\bf z}}\}$ be a generating set of $\mathcal{C}_{A}$ and let $1 \notin {\rm supp}({\bf x}^{{\bf v}})$. In the following cases $C$ has Cohen-Macaulay tangent cone at the origin. (i) $1 \in {\rm supp}({\bf x}^{{\bf w}})$, $1 \in {\rm supp}({\bf x}^{{\bf z}})$ and \begin{enumerate} \item $a_{3} \leq {\rm deg}({\bf x}^{{\bf w}})$, \item for every binomial $f=M-N \in \mathcal{I}$ with $1 \in {\rm supp}(M)$ we have that ${\rm deg}(N) \leq {\rm deg}(M)$ and \item for every monomial $M=x_{2}^{d_2}x_{3}^{d_3}x_{4}^{d_4}$, where $d_{2} \geq a_{2}$, $d_{3}<a_{3}$ and $d_{4}<a_{4}$, with $d_{2}n_{2}+d_{3}n_{3}+d_{4}n_{4} \in n_{1}+\mathcal{S}$ there exists a monomial $N$ with $1 \in {\rm supp}(N)$ such that $M-N \in I(C)$ and also ${\rm deg}(M) \leq {\rm deg}(N)$.
\end{enumerate}
(ii) $1 \in {\rm supp}({\bf x}^{{\bf w}})$, $1 \notin {\rm supp}({\bf x}^{{\bf z}})$ and \begin{enumerate} \item $a_{3} \leq {\rm deg}({\bf x}^{{\bf w}})$ and \item for every monomial $M=x_{2}^{d_2}x_{3}^{d_3}x_{4}^{d_4}$, where $d_{3}<a_{3}$, with $d_{2}n_{2}+d_{3}n_{3}+d_{4}n_{4} \in n_{1}+\mathcal{S}$ there exists a monomial $N$ with $1 \in {\rm supp}(N)$ such that $M-N \in I(C)$ and also ${\rm deg}(M) \leq {\rm deg}(N)$.
\end{enumerate}
(iii) $1 \notin {\rm supp}({\bf x}^{{\bf w}})$, $1 \in {\rm supp}({\bf x}^{{\bf z}})$ and for every monomial $M=x_{2}^{d_2}x_{3}^{d_3}x_{4}^{d_4}$, where $d_{4}<a_{4}$, with $d_{2}n_{2}+d_{3}n_{3}+d_{4}n_{4} \in n_{1}+\mathcal{S}$ there exists a monomial $N$ with $1 \in {\rm supp}(N)$ such that $M-N \in I(C)$ and also ${\rm deg}(M) \leq {\rm deg}(N)$.
\end{proposition}

\noindent \textbf{Proof.} (i) Let $M=x_{2}^{d_2}x_{3}^{d_3}x_{4}^{d_4}$, where $d_{2} \geq 0$, $d_{3} \geq 0$ and $d_{4} \geq 0$, with $d_{2}n_{2}+d_{3}n_{3}+d_{4}n_{4} \in n_{1}+\mathcal{S}$. Thus there exists a monomial $P$ such that $1 \in {\rm supp}(P)$ and $M-P \in I(C)$. Let $d_{3} \geq a_{3}$, then we consider the monomial $P=x_{2}^{d_{2}}x_{3}^{d_{3}-a_{3}}x_{4}^{d_4}{\bf x}^{{\bf w}}$. We have that $M-P \in I(C)$ and also $${\rm deg}(M)=d_{2}+d_{3}+d_{4} \leq d_{2}+d_{3}+d_{4}+({\rm deg}({\bf x}^{{\bf w}})-a_{3})={\rm deg}(P).$$ Let $d_{4} \geq a_{4}$, then we consider the monomial $P=x_{2}^{d_{2}}x_{3}^{d_{3}}x_{4}^{d_{4}-a_{4}}{\bf x}^{{\bf z}}$. We have that $M-P \in I(C)$ and also $${\rm deg}(M)=d_{2}+d_{3}+d_{4}<d_{2}+d_{3}+d_{4}+({\rm deg}({\bf x}^{{\bf z}})-a_{4})={\rm deg}(P).$$ Suppose that $d_{3}<a_{3}$ and $d_{4}<a_{4}$. If $d_{2} \geq a_{2}$, then from (3) we are done. Assume that $d_{2}<a_{2}$. There are 2 cases: (1) there exists a binomial $f=G-H \in \mathcal{I}$ with $1 \in {\rm supp}(G)$ such that $H$ divides $M$, so $M=H H'$. Note that ${\rm deg}(H) \geq {\rm deg}(G)$ from condition (2). Then $GH'-M \in I(C)$ and also ${\rm deg}(M) \leq {\rm deg}(G)+{\rm deg}(H')$.\\
(2) there exists no binomial $f=G-H$ in $\mathcal{I}$ such that $H$ divides $M$. Recall that $M-P \in I(C)$, $S \cup \mathcal{I}$ generates $I(C)$ and also $1 \in {\rm supp}({\bf x}^{{\bf w}})$ and $1 \in {\rm supp}({\bf x}^{{\bf z}})$. Then necessarily ${\bf x}^{{\bf u}}$ or/and ${\bf x}^{{\bf v}}$ divides $M$. Let us suppose that ${\bf x}^{{\bf u}}$ divides $M$, so $M={\bf x}^{{\bf u}}M'$. Let $P=M'x_{1}^{a_{1}}$, then $M-P \in I(C)$ and also $${\rm deg}(M)={\rm deg}({\bf x}^{{\bf u}})+{\rm deg}(M')<{\rm deg}(x_{1}^{a_1})+{\rm deg}(M')={\rm deg}(P).$$ Suppose now that ${\bf x}^{{\bf v}}$ divides $M$, so $M={\bf x}^{{\bf v}}M'$. Then the binomial $x_{2}^{a_{2}}M'-M$ belongs to $I(C)$ and also ${\rm deg}_{\mathcal{S}}(x_{2}^{a_{2}}M') \in n_{1}+\mathcal{S}$. Thus there exists a monomial $N$ such that $1 \in {\rm supp}(N)$, ${\rm deg}_{\mathcal{S}}(N)={\rm deg}_{\mathcal{S}}(x_{2}^{a_{2}}M')$ and also ${\rm deg}(x_{2}^{a_{2}}M') \leq {\rm deg}(N)$. Consequently $${\rm deg}(M)={\rm deg}({\bf x}^{{\bf v}}M')<{\rm deg}(x_{2}^{a_{2}}M') \leq {\rm deg}(N).$$

(ii) Let $M=x_{2}^{d_2}x_{3}^{d_3}x_{4}^{d_4}$, where $d_{2} \geq 0$, $d_{3} \geq 0$ and $d_{4} \geq 0$, with $d_{2}n_{2}+d_{3}n_{3}+d_{4}n_{4} \in n_{1}+\mathcal{S}$. Let $d_{3} \geq a_{3}$, then we consider the monomial $P=x_{2}^{d_{2}}x_{3}^{d_{3}-a_{3}}x_{4}^{d_4}{\bf x}^{{\bf w}}$. We have that $M-P \in I(C)$ and also $${\rm deg}(M)=d_{2}+d_{3}+d_{4} \leq d_{2}+d_{3}+d_{4}+({\rm deg}({\bf x}^{{\bf w}})-a_{3})={\rm deg}(P).$$ Let $d_{3}< a_{3}$, then from condition (2) we are done.\\
(iii) Let $M=x_{2}^{d_2}x_{3}^{d_3}x_{4}^{d_4}$, where $d_{2} \geq 0$, $d_{3} \geq 0$ and $d_{4} \geq 0$, with $d_{2}n_{2}+d_{3}n_{3}+d_{4}n_{4} \in n_{1}+\mathcal{S}$. Let $d_{4} \geq a_{4}$, then we consider the monomial $P=x_{2}^{d_{2}}x_{3}^{d_{3}}x_{4}^{d_{4}-a_{4}}{\bf x}^{{\bf z}}$. We have that $M-P \in I(C)$ and also $${\rm deg}(M)=d_{2}+d_{3}+d_{4}<d_{2}+d_{3}+d_{4}+({\rm deg}({\bf x}^{{\bf z}})-a_{4})={\rm deg}(P).$$ Suppose that $d_{4}<a_{4}$, then from the assumption we are done.

\section{The Gorenstein case}

In this section we will study the case that $C$ is a non-complete intersection Gorenstein monomial curve, i.e. the semigroup $\mathcal{S}=\{g_{1}n_{1}+\cdots+g_{4}n_{4}|g_{i} \in \mathbb{N}\}$ is symmetric.

\begin{theorem} (\cite{Bresinsky75}) Let $C$ be a monomial curve having the
parametrization $$x_1 = t^{n_1}, x_2 = t^{n_2}, x_3 = t^{n_3}, x_4 = t^{n_4}.$$
The semigroup $\mathcal{S}$ is symmetric and $C$ is a non-complete intersection curve if and only if $I(C)$ is minimally generated by the set
$$G = \{f_{1} = x_1^{a_1}- x_3^{a_{13}} x_4^{a_{14}}, f_2 = x_{2}^{a_2}- x_{1}^{a_{21}}x_{4}^{a_{24}}, f_3 = x_3^{a_{3}}-x_{1}^{a_{31}}x_{2}^{a_{32}},$$
$$f_4 = x_{4}^{a_4}-x_{2}^{a_{42}}x_{3}^{a_{43}}, f_5 = x_{3}^{a_{43}}
x_{1}^{a_{21}}-x_{2}^{a_{32}}x_4^{a_{14}}\}$$
where the polynomials $f_i$ are unique up to isomorphism and $0 < a_{ij} < a_{j}$.

\end{theorem}

\begin{remark} {\rm Bresinsky \cite{Bresinsky75} showed that $\mathcal{S}$ is symmetric and $I(C)$ is as in the previous theorem if and only if $n_1 =a_{2}a_{3}a_{14}+a_{32}a_{13}a_{24}$, $n_2 =a_{3}a_{4}a_{21} +a_{31}a_{43}a_{24}$, $n_3 =a_{1}a_{4}a_{32}+
a_{14}a_{42}a_{31}$, $n_4 =a_{1}a_{2}a_{43}+a_{42}a_{21}a_{13}$ with ${\rm gcd}(n_1, n_2, n_3, n_4) = 1$, $a_i > 1, 0 <a_{ij} < a_{j}$ for $1 \leq i \leq 4$ and $a_{1} =a_{21}+a_{31}$, $a_{2}= a_{32}+a_{42}$, $a_{3}=a_{13}+a_{43}$, $a_{4} =a_{14}+a_{24}$.}

\end{remark}
\begin{remark} (\cite{ArMe}) \label{BasicGorenstein} {\rm The above theorem implies that for any non-complete intersection Gorenstein monomial curve with embedding dimension four, the variables can be renamed to
obtain generators exactly of the given form, and this means that there are six
isomorphic possible permutations which can be considered within three cases:
\begin{enumerate}
\item[(1)] $f_1 = (1,(3, 4))$
\begin{enumerate}
\item[(a)] $f_2 = (2,(1, 4))$, $f_3 = (3,(1, 2))$, $f_4 = (4,(2, 3))$, $f_5 = ((1, 3),(2, 4))$
\item[(b)] $f_2 = (2,(1, 3))$, $f_3 = (3,(2, 4))$, $f_4 = (4,(1, 2))$, $f_5 = ((1, 4),(2, 3))$
\end{enumerate}
\item[(2)] $f_1 = (1,(2, 3))$
\begin{enumerate}
\item[(a)] $f_2 = (2,(3, 4))$, $f_3 = (3,(1, 4))$, $f_4 = (4,(1, 2))$, $f_5 = ((2, 4),(1, 3))$
\item[(b)] $f_2 = (2,(1, 4))$, $f_3 = (3,(2, 4))$, $f_4 = (4,(1, 3))$, $f_5 = ((1, 2),(4, 3))$
\end{enumerate}

\item[(3)] $f_1 = (1,(2, 4))$
\begin{enumerate}
\item[(a)] $f_2 = (2,(1, 3))$, $f_3 = (3,(1, 4))$, $f_4 = (4,(2, 3))$, $f_5 = ((1, 2),(3, 4))$
\item[(b)] $f_2 = (2,(3, 4))$, $f_3 = (3,(1, 2))$, $f_4 = (4,(1, 3))$, $f_5 = ((2, 3),(1, 4))$
\end{enumerate}
\end{enumerate}
Here, the notation $f_i = (i,(j, k))$ and $f_5 = ((i, j),(k,l))$ denote the generators
$f_i = x_{i}^{a_i}-x_{j}^{a_{ij}}x_{k}^{a_{ik}}$ and $f_5 = x_{i}^{a_{ki}}x_{j}^{a_{lj}}-
x_{k}^{a_{jk}}x_{l}^{a_{il}}$. Thus, given a Gorenstein monomial curve $C$, if we have the extra condition $n_1 < n_2 < n_3 < n_4$, then the generator set of $I(C)$ is exactly given by one of these six permutations.}

\end{remark}

\begin{remark} {\rm By \cite[Corollary 3.13]{KaOj} the toric ideal $I(C)$ of any non-complete intersection Gorenstein monomial curve $C$ is generated by its indispensable binomials.}
\end{remark}

Let $lex-inf$ be the total order on the monomials of $K[x_{1},\ldots,x_{4}]$ which is defined as follows: $${\bf x}^{\bf u}>_{lex-inf} {\bf x}^{\bf v} \Leftrightarrow {\bf x}^{\bf u}<_{lex} {\bf x}^{\bf v}$$ where $lex$ order is the lexicographic order such that $x_{1}$ is the largest variable in $K[x_{1},\ldots,x_{4}]$ with respect to $<_{lex}$.

\begin{proposition} The set $G=\{f_{1}, f_2, f_3, f_4, f_5\}$ is the reduced Gr\"obner basis of $I(C)$ with respect to an appropriate $lex-inf$ order.
\end{proposition}

\noindent \textbf{Proof.} Suppose that $I(C)$ is given as in case 1(a). Then $$f_{1}=x_1^{a_1}-x_3^{a_{13}} x_4^{a_{14}}, f_2 = x_{2}^{a_2}- x_{1}^{a_{21}}x_{4}^{a_{24}}, f_3 = x_3^{a_{3}}-x_{1}^{a_{31}}x_{2}^{a_{32}},$$ $$f_4 = x_{4}^{a_4}-x_{2}^{a_{42}}x_{3}^{a_{43}}, f_5 =x_{1}^{a_{21}}x_{3}^{a_{43}}-x_{2}^{a_{32}}x_4^{a_{14}}.$$ With respect to $lex-inf$ such that $x_{1}>_{lex} x_{2} >_{lex} x_{3} >_{lex} x_{4}$ we have that $lm(f_{1})=x_3^{a_{13}} x_4^{a_{14}}$, $lm(f_{2})=x_2^{a_{2}}$, $lm(f_{3})=x_3^{a_{3}}$, $lm(f_{4})=x_4^{a_{4}}$ and $lm(f_{5})=x_{2}^{a_{32}}x_4^{a_{14}}$. We will prove that $S(f_{i},f_{j}) \stackrel{G}  {\longrightarrow} 0$ for any pair $\{f_{i},f_{j}\}$. Since $lm(f_{1})$ and $lm(f_{2})$ are relatively prime, we get that $S(f_{1},f_{2}) \stackrel{G}  {\longrightarrow} 0$. Similarly $S(f_{2},f_{3}) \stackrel{G}  {\longrightarrow} 0$, $S(f_{2},f_{4}) \stackrel{G}  {\longrightarrow} 0$, $S(f_{3},f_{4}) \stackrel{G}  {\longrightarrow} 0$ and $S(f_{3},f_{5}) \stackrel{G}  {\longrightarrow} 0$. We have that $$S(f_{1},f_{3})=x_{1}^{a_{31}}x_{2}^{a_{32}}x_{4}^{a_{14}}-x_{1}^{a_{1}}x_{3}^{a_{43}} \stackrel{f_5}  {\longrightarrow} x_{1}^{a_{1}}x_{3}^{a_{43}}-x_{1}^{a_{1}}x_{3}^{a_{43}}=0$$ $$S(f_{1},f_{4})=x_{2}^{a_{42}}x_{3}^{a_{3}}-x_{1}^{a_{1}}x_{4}^{a_{24}} \stackrel{f_3}  {\longrightarrow} x_{1}^{a_{31}}x_{2}^{a_{2}}-x_{1}^{a_{1}}x_{4}^{a_{24}}  \stackrel{f_2}  {\longrightarrow} x_{1}^{a_{1}}x_{4}^{a_{24}}-x_{1}^{a_{1}}x_{4}^{a_{24}}=0$$ $$S(f_{1},f_{5})=x_{1}^{a_{21}}x_{3}^{a_{3}}-x_{1}^{a_{1}}x_{2}^{a_{32}} \stackrel{f_3}  {\longrightarrow} x_{1}^{a_{1}}x_{2}^{a_{32}}-x_{1}^{a_{1}}x_{2}^{a_{32}}=0$$ $$S(f_{2},f_{5})=x_{1}^{a_{21}}x_{2}^{a_{42}}x_{3}^{a_{43}}-x_{1}^{a_{21}}x_{4}^{a_{4}} \stackrel{f_4}  {\longrightarrow} x_{1}^{a_{21}}x_{2}^{a_{42}}x_{3}^{a_{43}}-x_{1}^{a_{21}}x_{2}^{a_{42}}x_{3}^{a_{43}}=0$$ $$S(f_{4},f_{5})=x_{1}^{a_{21}}x_{3}^{a_{43}}x_{4}^{a_{24}}-x_{2}^{a_{2}}x_{3}^{a_{43}} \stackrel{f_2}  {\longrightarrow} x_{1}^{a_{21}}x_{3}^{a_{43}}x_{4}^{a_{24}}-x_{1}^{a_{21}}x_{3}^{a_{43}}x_{4}^{a_{24}}=0.$$ Thus $G$ is a Gr\"obner basis for $I(C)$ with respect to $lex-inf$ such that $x_{1}>_{lex} x_{2} >_{lex} x_{3} >_{lex} x_{4}$. It is not hard to show that $G$ is a Gr\"obner basis for $I(C)$ with respect to $lex-inf$ such that 
\begin{enumerate} \item $x_{1}>_{lex} x_{2} >_{lex} x_{3} >_{lex} x_{4}$ in case 1(b). \item $x_{1}>_{lex} x_{3} >_{lex} x_{2} >_{lex} x_{4}$ in case 2(a). \item $x_{1}>_{lex} x_{2} >_{lex} x_{3} >_{lex} x_{4}$ in case 2(b). \item $x_{1}>_{lex} x_{2} >_{lex} x_{3} >_{lex} x_{4}$ in case 3(a). \item $x_{1}>_{lex} x_{3} >_{lex} x_{2} >_{lex} x_{4}$ in case 3(b).

\end{enumerate}

The Apery set $Q$ of the semigroup $\mathcal{S}$ relative to $\{n_{1}\}$ is defined by $$Q=\{q \in \mathcal{S}|q-n_{1} \notin \mathcal{S}\}.$$ Using Lemma 1.2 in \cite{Cas} we get the following.

\begin{corollary} Let $B$ be the set of monomials $x_{2}^{u_2}x_{3}^{u_3}x_{4}^{u_4}$ in the polynomial ring $K[x_{2},x_{3},x_{4}]$ which are not divisible by any of the monomials of the set \begin{enumerate} \item $\{x_3^{a_{13}} x_4^{a_{14}},x_{2}^{a_2},x_3^{a_{3}},x_{4}^{a_4},x_{2}^{a_{32}}x_4^{a_{14}}\}$ in case 1(a). \item $\{x_3^{a_{13}} x_4^{a_{14}},x_{2}^{a_2},x_3^{a_{3}},x_{4}^{a_4},x_{2}^{a_{42}}x_3^{a_{13}}\}$ in case 1(b). \item $\{x_2^{a_{12}} x_3^{a_{13}},x_{2}^{a_2},x_3^{a_{3}},x_{4}^{a_4},x_{2}^{a_{12}}x_4^{a_{34}}\}$ in case 2(a). \item $\{x_2^{a_{12}} x_3^{a_{13}},x_{2}^{a_2},x_3^{a_{3}},x_{4}^{a_4},x_{3}^{a_{13}}x_4^{a_{24}}\}$ in case 2(b). \item $\{x_2^{a_{12}} x_4^{a_{14}},x_{2}^{a_2},x_3^{a_{3}},x_{4}^{a_4}, x_{3}^{a_{23}}x_4^{a_{14}}\}$ in case 3(a). \item $\{x_2^{a_{12}} x_4^{a_{14}},x_{2}^{a_2},x_3^{a_{3}},x_{4}^{a_4},x_{2}^{a_{12}}x_3^{a_{43}}\}$ in case 3(b).
\end{enumerate}
Then $$Q=\{m \in \mathcal{S}|m=\sum_{i=2}^{4}u_{i}n_{i} \ \textrm{where} \ x_{2}^{u_2}x_{3}^{u_3}x_{4}^{u_4} \in B\}.$$
\end{corollary}

\begin{theorem} Suppose that $I(C)$ is given as in case 1(a). Then $C$ has Cohen-Macaulay tangent cone at the origin if and only if $a_{2} \leq a_{21}+a_{24}$.

\end{theorem}

\noindent \textbf{Proof.} In this case $I(C)$ is minimally generated by the set $$G=\{f_{1} = x_1^{a_1}- x_3^{a_{13}} x_4^{a_{14}}, f_2 = x_{2}^{a_2}- x_{1}^{a_{21}}x_{4}^{a_{24}}, f_3 = x_3^{a_{3}}-x_{1}^{a_{31}}x_{2}^{a_{32}},$$ $$f_4 = x_{4}^{a_4}-x_{2}^{a_{42}}x_{3}^{a_{43}}, f_5 = x_{1}^{a_{21}}x_{3}^{a_{43}}-x_{2}^{a_{32}}x_4^{a_{14}}\}.$$ If $a_{2} \leq a_{21}+a_{24}$, then we have, from Theorem 2.8 in \cite{ArMe}, that the curve $C$ has Cohen-Macaulay tangent cone at the origin. Conversely suppose that $C$ has Cohen-Macaulay tangent cone at the origin. Since $I(C)$ is generated by the indispensable binomials, every binomial $f_i$, $1 \leq i \leq 5$, is indispensable of $I(C)$. In particular the binomial $f_2$ is indispensable of $I(C)$. If there exists a monomial $N \neq x_{1}^{a_{21}}x_{4}^{a_{24}}$ such that $g=x_{2}^{a_2}-N$ belongs to $I(C)$, then we can replace $f_{2}$ in $G$ by the binomials $g$ and $N-x_{1}^{a_{21}}x_{4}^{a_{24}} \in I(C)$, a contradiction to the fact that $f_2$ is indispensable. Thus $N=x_{1}^{a_{21}}x_{4}^{a_{24}}$. But $a_{2}n_{2} \in n_{1}+\mathcal{S}$ and therefore we have, from Theorem \ref{BasicHer}, that $a_{2} \leq a_{21}+a_{24}$.

\begin{remark} {\rm Suppose that $I(C)$ is given as in case 1(b). (1) It holds that $a_{1}>a_{13}+a_{14}$ and $a_{4}<a_{41}+a_{42}$.\\ (2) If $a_{42} \leq a_{32}$, then $x_{3}^{a_{3}+a_{13}}-x_{1}^{a_{21}}x_{2}^{a_{32}-a_{42}}x_{4}^{2a_{34}} \in I(C)$.\\(3) If $a_{14} \leq a_{34}$, then the binomial $x_{3}^{a_{3}+a_{13}}-x_{1}^{a_{1}}x_{2}^{a_{32}}x_{4}^{a_{34}-a_{14}}$ belongs to $I(C)$.}
\end{remark}

\begin{proposition} \label{Prop1(b)} Suppose that $I(C)$ is given as in case 1(b). Then $C$ has Cohen-Macaulay tangent cone at the origin if and only if \begin{enumerate} \item $a_{2} \leq a_{21}+a_{23}$, \item $a_{42}+a_{13} \leq a_{21} +a_{34}$ and \item for every monomial $M=x_{2}^{u_{2}}x_{3}^{u_{^3}}x_{4}^{u_{4}}$, where $u_{2}<a_{42}$, $u_{3} \geq a_{3}$ and $u_{4}<a_{14}$, with $u_{2}n_{2}+u_{3}n_{3}+u_{4}n_{4} \in n_{1}+\mathcal{S}$ there exists a monomial $N$ with $1 \in {\rm supp}(N)$ such that $M-N \in I(C)$ and also ${\rm deg}(M) \leq {\rm deg}(N)$.
\end{enumerate}

\end{proposition}

\noindent \textbf{Proof.} In this case $I(C)$ is minimally generated by the set $$G=\{f_{1} = x_1^{a_1}- x_3^{a_{13}} x_4^{a_{14}}, f_2 = x_{2}^{a_2}- x_{1}^{a_{21}}x_{3}^{a_{23}}, f_3 = x_3^{a_{3}}-x_{2}^{a_{32}}x_{4}^{a_{34}},$$ $$f_4 = x_{4}^{a_4}-x_{1}^{a_{41}}x_{2}^{a_{42}}, f_5 = x_{1}^{a_{21}}x_{4}^{a_{34}}-x_{2}^{a_{42}}x_3^{a_{13}}\}.$$
Suppose that $C$ has Cohen-Macaulay tangent cone at the origin. Since $I(C)$ is generated by the indispensable binomials, every binomial $f_i$, $1 \leq i \leq 5$, is indispensable of $I(C)$. In particular the binomials $f_2$ and $f_5$ are indispensable of $I(C)$. Therefore both inequalities $a_{2} \leq a_{21}+a_{23}$ and $a_{42}+a_{13} \leq a_{21} +a_{34}$ hold. By Theorem \ref{BasicHer} condition (3) is also true.

Conversely, from Proposition \ref{1(b)} (iii), it is enough to consider a monomial $M=x_{2}^{u_2}x_{3}^{u_3}x_{4}^{u_4}$, where $u_2<a_{2}$, $u_3 \geq 0$ and $u_4<a_{4}$, with the property: there exists at least one monomial $P$ such that $1 \in {\rm supp}(P)$ and also $M-P$ is in $I(C)$. Suppose that $u_{3} \geq a_{3}$. If $u_{4} \geq a_{14}$, then we let $P=x_{1}^{a_{1}}x_{2}^{u_{2}}x_{3}^{u_{3}-a_{13}}x_{4}^{u_{4}-a_{14}}$. So we have that $M-P \in I(C)$ and also ${\rm deg}(M)<{\rm deg}(P)$ since $a_{13}+a_{14}<a_{1}$. Similarly if $u_{2} \geq a_{42}$, then we let $P=x_{1}^{a_{21}}x_{2}^{u_{2}-a_{42}}x_{3}^{u_{3}-a_{13}}x_{4}^{u_{4}+a_{34}}$. So we have that $M-P \in I(C)$ and also ${\rm deg}(M) \leq {\rm deg}(P)$. If both inequalities $u_{4}<a_{14}$ and $u_{2}<a_{42}$ hold, then condition (3) implies that there exists a monomial $N$ with $1 \in {\rm supp}(N)$ such that $M-N \in I(C)$ and also ${\rm deg}(M) \leq {\rm deg}(N)$.

Suppose now that $u_3<a_3$. Recall that $M-P \in I(C)$ and $G$ generates $I(C)$. Then $M$ is divided by at least one of the monomials $x_{2}^{a_{42}}x_3^{a_{13}}$, $x_{3}^{a_{13}}x_4^{a_{14}}$ and $x_{2}^{a_{32}}x_4^{a_{34}}$. If $M$ is divided by $x_{2}^{a_{42}}x_3^{a_{13}}$, then $M=x_2^{a_{42}+p}x_3^{a_{13}+q}x_{4}^{u_4}$, for some non-negative integers $p$ and $q$, so $M-x_1^{a_{21}}x_2^{p}x_3^{q}x_{4}^{a_{34}+u_{4}} \in I(C)$ and also ${\rm deg}(M) \leq {\rm deg}(x_1^{a_{21}}x_2^{p}x_3^{q}x_{4}^{a_{34}+u_{4}})$. If $M$ is divided by $x_{3}^{a_{13}}x_4^{a_{14}}$, then $M=x_{2}^{u_2}x_3^{a_{13}+p}x_4^{a_{14}+q}$, for some non-negative integers $p$ and $q$, and therefore the binomial $M-x_1^{a_{1}}x_{2}^{u_2}x_3^{p}x_4^{q} \in I(C)$ and also ${\rm deg}(M)<{\rm deg}(x_1^{a_{1}}x_{2}^{u_2}x_3^{p}x_4^{q})$. Assume that neither $x_{2}^{a_{32}}x_4^{a_{34}}$ nor $x_{3}^{a_{13}}x_4^{a_{14}}$ divides $M$. Then necessarily  $x_{2}^{a_{42}}x_3^{a_{13}}$ divides $M$. But $M$ is not divided by any leading monomial of $G$ with respect to $lex-inf$ such that $x_{1}>_{lex} x_{2} >_{lex} x_{3} >_{lex} x_{4}$. Thus $m=u_{2}n_{2}+u_{3}n_{3}+u_{4}n_{4}$ is in $Q$, a contradiction to the fact that $m-n_{1} \in \mathcal{S}$. Therefore, from Proposition \ref{1(b)}, $C$ has Cohen-Macaulay tangent cone at the origin.

\begin{proposition} \label{CasesCohen1(b)} Suppose that $I(C)$ is given as in case 1(b). Assume that $C$ has Cohen-Macaulay tangent cone at the origin and also $a_{42} \leq a_{32}$. \begin{enumerate} \item If $a_{34}<a_{14}$, then $a_{3}+a_{13} \leq a_{21}+a_{32}-a_{42}+2a_{34}$.
\item If $a_{14} \leq a_{34}$, then $a_{3}+a_{13} \leq a_{1}+a_{32}+a_{34}-a_{14}$.
\end{enumerate}

\end{proposition}

\noindent \textbf{Proof.} First of all we note that $a_{21}+a_{32}-a_{42}+2a_{34}<a_{1}+a_{32}+a_{34}-a_{14}$, since $a_{14}+a_{34}<a_{41}+a_{42}=a_{1}-a_{21}+a_{42}$, so $a_{34}<a_{1}-a_{21}+a_{42}-a_{14}$ and therefore $$a_{21}+a_{32}-a_{42}+2a_{34}=a_{34}+(a_{21}+a_{32}-a_{42}+a_{34})<a_{1}+a_{32}+a_{34}-a_{14}.$$ Since $C$ has Cohen-Macaulay tangent cone at the origin, there is a monomial $M$ with $1 \in {\rm supp}(M)$ such that $x_{3}^{a_{3}+a_{13}}-M \in I(C)$ and also $a_{3}+a_{13} \leq {\rm deg}(M)$. Now $x_{3}^{a_{3}+a_{13}}-M \in I(C)$, so $x_{3}^{a_{3}+a_{13}}-M=\sum_{i=1}^{5}H_{i}f_{i}$ for some $H_{i} \in K[x_{1},\ldots,x_{4}]$. Then $x_{3}^{a_{3}+a_{13}}$ arise as a term in the sum $\sum_{i=1}^{5}H_{i}f_{i}$, so $Q=-x_{2}^{a_{32}}x_{3}^{a_{13}}x_{4}^{a_{34}}$ is a term in the sum $\sum_{i=1}^{5}H_{i}f_{i}$, therefore $Q$ should be canceled with another term of the above sum. We distinguish the following cases:\\ (1) $a_{34}<a_{14}$. Suppose that $M \neq x_{1}^{a_{21}}x_{2}^{a_{32}-a_{42}}x_{4}^{2a_{34}}$. Note that $x_{2}^{a_{32}}x_{3}^{a_{13}}x_{4}^{a_{34}}$ is divided only by $x_{2}^{a_{42}}x_{3}^{a_{13}}$ and $x_{2}^{a_{32}}x_{4}^{a_{34}}$. Thus the term $T=-x_{1}^{a_{21}}x_{2}^{a_{32}-a_{42}}x_{4}^{2a_{34}}$ arise in the sum $\sum_{i=1}^{5}H_{i}f_{i}$ and it should be canceled with another term of the sum. Consequently the monomial $x_{1}^{a_{41}}x_{2}^{a_{42}}$ must divide $-T$, so $a_{41} \leq a_{21}$ and $a_{42} \leq a_{32}-a_{42}$. Therefore the term $R=-x_{1}^{a_{21}-a_{41}}x_{2}^{a_{32}-2a_{42}}x_{4}^{a_{4}+2a_{34}}$ arise in the sum $\sum_{i=1}^{5}H_{i}f_{i}$. Note that ${\rm deg}(-R)<a_{21}+a_{32}-a_{42}+2a_{34}$. If $M=-R$, then we are done. Otherwise $R$ should be canceled with another term of the above sum. Thus $x_{1}^{a_{41}}x_{2}^{a_{42}}$ must divide $-R$, so $a_{41} \leq a_{21}-a_{41}$ and $a_{42} \leq a_{32}-2a_{42}$. Consequently the term $V=-x_{1}^{a_{21}-2a_{41}}x_{2}^{a_{32}-3a_{42}}x_{4}^{2a_{4}+2a_{34}}$ arise in the sum $\sum_{i=1}^{5}H_{i}f_{i}$. If $M=-V$, then we are done. Otherwise $V$ should be canceled with another term of the above sum. Continuing this way we finally reach a contradiction. Thus $M=x_{1}^{a_{21}}x_{2}^{a_{32}-a_{42}}x_{4}^{2a_{34}}$.\\ (2) $a_{14} \leq a_{34}$. Suppose that $M \neq x_{1}^{a_{21}}x_{2}^{a_{32}-a_{42}}x_{4}^{2a_{34}}$ and also $M \neq x_{1}^{a_{1}}x_{2}^{a_{32}}x_{4}^{a_{34}-a_{14}}$. Note that $x_{2}^{a_{32}}x_{3}^{a_{13}}x_{4}^{a_{34}}$ is divided only by $x_{2}^{a_{32}}x_{4}^{a_{34}}$, $x_{3}^{a_{13}}x_{4}^{a_{14}}$ and $x_{2}^{a_{42}}x_{3}^{a_{13}}$. If the term $T=-x_{1}^{a_{1}}x_{2}^{a_{32}}x_{4}^{a_{34}-a_{14}}$ arise in the sum $\sum_{i=1}^{5}H_{i}f_{i}$, then it should be canceled with another term of the above sum. Note that $x_{1}^{a_{1}}x_{2}^{a_{32}}x_{4}^{a_{34}-a_{14}}$ is divided by $x_{1}^{a_{41}}x_{2}^{a_{42}}$, so the term $T=-x_{1}^{a_{21}}x_{2}^{a_{32}-a_{42}}x_{4}^{a_{4}+a_{34}-a_{14}}=-x_{1}^{a_{21}}x_{2}^{a_{32}-a_{42}}x_{4}^{2a_{34}}$ arise in the sum $\sum_{i=1}^{5}H_{i}f_{i}$. Since $M \neq -T$, $T$ should be canceled with another term of the above sum. Then $a_{41} \leq a_{21}$ and $a_{42} \leq a_{32}-a_{42}$. Thus the term $R=-x_{1}^{a_{21}-a_{41}}x_{2}^{a_{32}-2a_{42}}x_{4}^{a_{4}+2a_{34}}$ arise in the sum $\sum_{i=1}^{5}H_{i}f_{i}$. If $M=-R$, then we are done. Otherwise $R$ should be canceled with another term of the above sum. Continuing this way we finally reach a contradiction. Thus either $M =x_{1}^{a_{21}}x_{2}^{a_{32}-a_{42}}x_{4}^{2a_{34}}$ or $M=x_{1}^{a_{1}}x_{2}^{a_{32}}x_{4}^{a_{34}-a_{14}}$.

\begin{theorem} Suppose that $I(C)$ is given as in case 1(b) and also that $a_{42} \leq a_{32}$. Then $C$ has Cohen-Macaulay tangent cone at the origin if and only if \begin{enumerate} \item $a_{2} \leq a_{21}+a_{23}$, \item $a_{42}+a_{13} \leq a_{21} +a_{34}$
and \item either $a_{34}<a_{14}$ and $a_{3}+a_{13} \leq a_{21}+a_{32}-a_{42}+2a_{34}$ or $a_{14} \leq a_{34}$ and $a_{3}+a_{13} \leq a_{1}+a_{32}+a_{34}-a_{14}$.
\end{enumerate}

\end{theorem}

\noindent \textbf{Proof.} $(\Longrightarrow)$ From Proposition \ref{Prop1(b)} we have that conditions (1) and (2) are true. From Proposition \ref{CasesCohen1(b)} the condition (3) is also true.\\
$(\Longleftarrow)$ From Proposition \ref{Prop1(b)} it is enough to consider a monomial $N=x_{2}^{u_2}x_{3}^{u_3}x_{4}^{u_4}$, where $u_2<a_{42}$, $u_3 \geq a_{3}$ and $u_4<a_{14}$, with the property: there exists at least one monomial $P$ such that $1 \in {\rm supp}(P)$ and also $N-P$ is in $I(C)$. Suppose that $u_{3} \geq a_{3}+a_{13}$ and let $M$ denote either the monomial $x_{1}^{a_{21}}x_{2}^{a_{32}-a_{42}}x_{4}^{2a_{34}}$ when $a_{34}<a_{14}$ or the monomial $x_{1}^{a_{1}}x_{2}^{a_{32}}x_{4}^{a_{34}-a_{14}}$ when $a_{14} \leq a_{34}$. Let $P=x_{2}^{u_2}x_{3}^{u_{3}-a_{3}-a_{13}}x_{4}^{u_4}M$. We have that $N-P \in I(C)$ and $${\rm deg}(N)=u_{2}+u_{3}+u_{4} \leq u_{2}+u_{3}+u_{4}+{\rm deg}(M)-a_{3}-a_{13}={\rm deg}(P).$$ It suffices to consider the case that $u_{3}-a_{3}<a_{13}$. Recall that $G=\{f_{1},\ldots,f_{5}\}$ generates $I(C)$. The binomial $N-P$ belongs to $I(C)$, so $N-P=\sum_{i=1}^{5}H_{i}f_{i}$ for some polynomials $H_{i} \in K[x_{1},\ldots,x_{4}]$ and therefore $N$ is a term in the sum $\sum_{i=1}^{5}H_{i}f_{i}$. Note that $N$ is not divided by any of the monomials $x_3^{a_{13}} x_4^{a_{14}}$, $x_{2}^{a_2}$, $x_{2}^{a_{32}}x_{4}^{a_{34}}$, $x_{4}^{a_4}$ and $x_{2}^{a_{42}}x_3^{a_{13}}$. Now the monomial $N$ is divided by the monomial $x_{3}^{a_{3}}$, so $Q=-x_{2}^{u_{2}+a_{32}}x_{3}^{u_{3}-a_{3}}x_{4}^{u_{4}+a_{34}}$ is a term in the sum $\sum_{i=1}^{5}H_{i}f_{i}$ and should be canceled with another term of the above sum. Remark that $u_{2}+a_{32}<a_{2}$ and $u_{4}+a_{34}<a_{4}$. Thus $x_{2}^{a_{42}}x_{3}^{a_{13}}$ divides $-Q$, so $u_{3}-a_{3} \geq a_{13}$ a contradiction.

\begin{proposition} \label{Difficult1} Suppose that $I(C)$ is given as in case 1(b) and also that $a_{32}<a_{42}$. If $C$ has Cohen-Macaulay tangent cone at the origin, then \begin{enumerate} \item $a_{2} \leq a_{21}+a_{23}$, \item $a_{42}+a_{13} \leq a_{21} +a_{34}$ and \item $a_{3}+a_{13} \leq a_{1}+a_{32}+a_{34}-a_{14}$.
\end{enumerate}

\end{proposition}

\noindent \textbf{Proof.} From Proposition \ref{Prop1(b)} we have that conditions (1) and (2) are true. Suppose first that $a_{14} \leq a_{34}$. Again from Proposition \ref{Prop1(b)}, for $u_{2}=u_{4}=0$ and $u_{3}=a_{3}+a_{13}$, we deduce that there exists a monomial $M$ with $1 \in {\rm supp}(M)$ such that $x_{3}^{a_{3}+a_{13}}-M \in I(C)$ and $a_{3}+a_{13} \leq {\rm deg}(M)$. Since the binomial $x_{3}^{a_{3}+a_{13}}-M$ belongs to $I(C)$, we have that $x_{3}^{a_{3}+a_{13}}-M=\sum_{i=1}^{5}H_{i}f_{i}$ for some polynomials $H_{i} \in K[x_{1},\ldots,x_{4}]$. Then $x_{3}^{a_{3}+a_{13}}$ arise as a term in the sum $\sum_{i=1}^{5}H_{i}f_{i}$, so $Q=-x_{2}^{a_{32}}x_{3}^{a_{13}}x_{4}^{a_{34}}$ is a term in the sum $\sum_{i=1}^{5}H_{i}f_{i}$, therefore $Q$ should be canceled with another term of the above sum. Thus $T=-x_{1}^{a_{1}}x_{2}^{a_{32}}x_{4}^{a_{34}-a_{14}}$ is a term in the sum $\sum_{i=1}^{5}H_{i}f_{i}$. Suppose that $M \neq x_{1}^{a_{1}}x_{2}^{a_{32}}x_{4}^{a_{34}-a_{14}}$. Then $T$ should be canceled with another term of the above sum, a contradiction. Thus $M=x_{1}^{a_{1}}x_{2}^{a_{32}}x_{4}^{a_{34}-a_{14}}$ and also $a_{3}+a_{13} \leq a_{1}+a_{32}+a_{34}-a_{14}$. Suppose now that $a_{14}>a_{34}$. Note that $x_{1}^{a_{1}}x_{2}^{a_{32}}-x_{3}^{a_{3}+a_{13}}x_{4}^{a_{14}-a_{34}}\in I(C)$. Thus $a_{1}+a_{32}>a_{3}+a_{13}+a_{14}-a_{34}$.
\begin{theorem} Suppose that $I(C)$ is given as in case 1(b) and also that $a_{32}<a_{42}$. Assume that $a_{14} \leq a_{34}$. Then $C$ has Cohen-Macaulay tangent cone at the origin if and only if \begin{enumerate} \item $a_{2} \leq a_{21}+a_{23}$, \item $a_{42}+a_{13} \leq a_{21} +a_{34}$ and \item $a_{3}+a_{13} \leq a_{1}+a_{32}+a_{34}-a_{14}$.
\end{enumerate}

\end{theorem}

\noindent \textbf{Proof.} $(\Longrightarrow)$ From Proposition \ref{Difficult1} we have that conditions (1), (2) and (3) are true.\\
$(\Longleftarrow)$ From Proposition \ref{Prop1(b)} it is enough to consider a monomial $N=x_{2}^{u_2}x_{3}^{u_3}x_{4}^{u_4}$, where $u_2<a_{42}$, $u_3 \geq a_{3}$ and $u_4<a_{14}$, with the property: there exists at least one monomial $P$ such that $1 \in {\rm supp}(P)$ and also $N-P$ is in $I(C)$. Suppose that $u_3 \geq a_{3}+a_{13}$. Let $P=x_{1}^{a_{1}}x_{2}^{u_{2}+a_{32}}x_{3}^{u_{3}-a_{3}-a_{13}}x_{4}^{u_{4}+a_{34}-a_{14}}$. We have that $N-P \in I(C)$ and $${\rm deg}(N)=u_{2}+u_{3}+u_{4} \leq u_{2}+u_{3}+u_{4}+a_{1}+a_{32}+a_{34}-a_{14}-a_{3}-a_{13}={\rm deg}(P).$$ It suffices to assume that $u_3-a_{3}<a_{13}$. Since the binomial $N-P$ belongs to $I(C)$, we have that $N-P=\sum_{i=1}^{5}H_{i}f_{i}$ for some polynomials $H_{i} \in K[x_{1},\ldots,x_{4}]$. So $N$ is a term in the sum $\sum_{i=1}^{5}H_{i}f_{i}$. Then $T=-x_{2}^{u_{2}+a_{32}}x_{3}^{u_{3}-a_{3}}x_{4}^{u_{4}+a_{34}}$ is a term in the sum $\sum_{i=1}^{5}H_{i}f_{i}$ and it should be canceled with another term of the above sum. Remark that $u_{2}+a_{32}<a_{2}$ and $u_{4}+a_{34}<a_{4}$. Thus $x_{2}^{a_{42}}x_{3}^{a_{13}}$ divides $-T$, so $u_{3}-a_{3} \geq a_{13}$ a contradiction.

\begin{example} {\rm Consider $n_{1}=1199$, $n_{2}=2051$, $n_{3}=2352$ and $n_{4}=3032$. The toric ideal $I(C)$ is minimally generated by the set $$G=\{x_1^{16}- x_{3}^{3}x_{4}^{4}, x_{2}^{19}- x_{1}^{7}x_{3}^{13}, x_3^{16}-x_{2}^{8}x_{4}^{7}, x_{4}^{11}-x_{1}^{9}x_{2}^{11}, x_{1}^{7}x_{4}^{7}-x_{2}^{11}x_{3}^{3}\}.$$ Here $a_{1}=16$, $a_{32}=8$, $a_{42}=11$, $a_{14}=4$, $a_{34}=7$, $a_{13}=3$ and $a_{3}=16$. Note that $a_{2}=19<20=a_{21}+a_{23}$ and $a_{42}+a_{13}=14=a_{21}+a_{34}$. We have that $a_{3}+a_{13}=19<27=a_{1}+a_{32}+a_{34}-a_{14}$. Thus $C$ has Cohen-Macaulay tangent cone at the origin.}

\end{example}

\begin{remark} {\rm Suppose that $I(C)$ is given as in case 2(a).\\ (1) It holds that $a_{1}>a_{12}+a_{13}$, $a_{2}>a_{23}+a_{24}$ and $a_{4}<a_{41}+a_{42}$.\\(2) If $a_{34} \leq a_{24}$, then $x_{2}^{a_{2}+a_{12}}-x_{1}^{a_{41}}x_{3}^{2a_{23}}x_{4}^{a_{24}-a_{34}} \in I(C)$.\\(3) If $a_{13}\leq a_{23}$, then the binomial $x_{2}^{a_{2}+a_{12}}-x_{1}^{a_{1}}x_{3}^{a_{23}-a_{13}}x_{4}^{a_{24}}$ belongs to $I(C)$.}
\end{remark}

\begin{proposition} \label{prop2(a)} Suppose that $I(C)$ is given as in case 2(a). Then $C$ has Cohen-Macaulay tangent cone at the origin if and only if \begin{enumerate} \item $a_{3} \leq a_{31}+a_{34}$, \item $a_{12}+a_{34} \leq a_{41} +a_{23}$ and \item for every monomial $M=x_{2}^{u_{2}}x_{3}^{u_{^3}}x_{4}^{u_{4}}$, where $u_{2} \geq a_{2}$, $u_{3}<a_{13}$ and $u_{4}<a_{34}$, with $u_{2}n_{2}+u_{3}n_{3}+u_{4}n_{4} \in n_{1}+\mathcal{S}$ there exists a monomial $N$ with $1 \in {\rm supp}(N)$ such that $M-N \in I(C)$ and also ${\rm deg}(M) \leq {\rm deg}(N)$.
\end{enumerate}

\end{proposition}

\noindent \textbf{Proof.} In this case $I(C)$ is minimally generated by the set $$G=\{f_{1} = x_1^{a_1}- x_2^{a_{12}} x_3^{a_{13}}, f_2 = x_{2}^{a_2}- x_{3}^{a_{23}}x_{4}^{a_{24}}, f_3 = x_3^{a_{3}}-x_{1}^{a_{31}}x_{4}^{a_{34}},$$ $$f_4 = x_{4}^{a_4}-x_{1}^{a_{41}}x_{2}^{a_{42}}, f_5 = x_{1}^{a_{41}}x_{3}^{a_{23}}-x_{2}^{a_{12}}x_4^{a_{34}}\}.$$
Suppose that $C$ has Cohen-Macaulay tangent cone at the origin. Since $I(C)$ is generated by the indispensable binomials, every binomial $f_i$, $1 \leq i \leq 5$, is indispensable of $I(C)$. In particular the binomials $f_3$ and $f_5$ are indispensable of $I(C)$. Therefore the inequalities $a_{3} \leq a_{31}+a_{34}$ and $a_{12}+a_{34} \leq a_{41} +a_{23}$ hold. By Theorem \ref{BasicHer} condition (3) is also true.

To prove the converse statement, from Proposition \ref{2(a)} (i), it is enough to consider a monomial $M=x_{2}^{u_2}x_{3}^{u_3}x_{4}^{u_4}$, where $u_2 \geq a_{2}$, $u_3<a_{3}$ and $u_4<a_{4}$, with the property: there exists at least one monomial $P$ such that $1 \in {\rm supp}(P)$ and also $M-P$ is in $I(C)$. If $u_{3} \geq a_{13}$, then we let $P=x_{1}^{a_{1}}x_{2}^{u_{2}-a_{12}}x_{3}^{u_{3}-a_{13}}x_{4}^{u_{4}}$. So we have that $M-P \in I(C)$ and also ${\rm deg}(M)<{\rm deg}(P)$. Similarly if $u_{4} \geq a_{34}$, then we let $P=x_{1}^{a_{41}}x_{2}^{u_{2}-a_{12}}x_{3}^{u_{3}+a_{23}}x_{4}^{u_{4}-a_{34}}$. So we have that $M-P \in I(C)$ and also ${\rm deg}(M) \leq {\rm deg}(P)$. If both conditions $u_{3}<a_{13}$ and $u_{4}<a_{34}$ hold, then condition (3) implies that there exists a monomial $N$ with $1 \in {\rm supp}(N)$ such that $M-N \in I(C)$ and also ${\rm deg}(M) \leq {\rm deg}(N)$. Therefore, from Proposition \ref{2(a)}, $C$ has Cohen-Macaulay tangent cone at the origin.

The proof of the next proposition is similar to that of Proposition \ref{CasesCohen1(b)} and therefore it is omitted.

\begin{proposition} \label{CasesCohen2(a)} Suppose that $I(C)$ is given as in case 2(a). Assume that $C$ has Cohen-Macaulay tangent cone at the origin and also $a_{34} \leq a_{24}$. \begin{enumerate} \item If $a_{23}<a_{13}$, then $a_{2}+a_{12} \leq a_{41}+2a_{23}+a_{24}-a_{34}$.
\item If $a_{13} \leq a_{23}$, then $a_{2}+a_{12} \leq a_{1}+a_{23}-a_{13}+a_{24}$.
\end{enumerate}

\end{proposition}

\begin{theorem} Suppose that $I(C)$ is given as in case 2(a) and also that $a_{34} \leq a_{24}$. Then $C$ has Cohen-Macaulay tangent cone at the origin if and only if \begin{enumerate} \item $a_{3} \leq a_{31}+a_{34}$, \item $a_{12}+a_{34} \leq a_{41} +a_{23}$
and \item either $a_{23}<a_{13}$ and $a_{2}+a_{12} \leq a_{41}+2a_{23}+a_{24}-a_{34}$ or $a_{13} \leq a_{23}$ and $a_{2}+a_{12} \leq a_{1}+a_{23}-a_{13}+a_{24}$.
\end{enumerate}

\end{theorem}

\noindent \textbf{Proof.} $(\Longrightarrow)$ From Proposition \ref{prop2(a)} we have that conditions (1) and (2) are true. From Proposition \ref{CasesCohen2(a)} the condition (3) is also true.\\
$(\Longleftarrow)$ From Proposition \ref{prop2(a)} it is enough to consider a monomial $N=x_{2}^{u_2}x_{3}^{u_3}x_{4}^{u_4}$, where $u_2 \geq a_{2}$, $u_3<a_{13}$ and $u_4<a_{34}$, with the property: there exists at least one monomial $P$ such that $1 \in {\rm supp}(P)$ and also $N-P$ is in $I(C)$. Suppose that $u_{2} \geq a_{2}+a_{12}$ and let $M$ denote either the monomial $x_{1}^{a_{41}}x_{3}^{2a_{23}}x_{4}^{a_{24}-a_{34}}$ when $a_{23}<a_{13}$ or the monomial $x_{1}^{a_{1}}x_{3}^{a_{23}-a_{13}}x_{4}^{a_{24}}$ when $a_{13} \leq a_{23}$. Let $P=x_{2}^{u_{2}-a_{2}-a_{12}}x_{3}^{u_{3}}x_{4}^{u_4}M$. We have that $N-P \in I(C)$ and $${\rm deg}(N)=u_{2}+u_{3}+u_{4} \leq u_{2}+u_{3}+u_{4}+{\rm deg}(M)-a_{2}-a_{12}={\rm deg}(P).$$ It suffices to consider the case that $u_{2}-a_{2}<a_{12}$. Since the binomial $N-P$ belongs to $I(C)$, we have that $N-P=\sum_{i=1}^{5}H_{i}f_{i}$ for some polynomials $H_{i} \in K[x_{1},\ldots,x_{4}]$. Now the monomial $N$ is divided by the monomial $x_{2}^{a_{2}}$, so $Q=-x_{2}^{u_{2}-a_{2}}x_{3}^{u_{3}+a_{23}}x_{4}^{u_{4}+a_{24}}$ is a term in the sum $\sum_{i=1}^{5}H_{i}f_{i}$ and should be canceled with another term of the above sum. Remark that $u_{3}+a_{23}<a_{3}$ and $u_{4}+a_{24}<a_{4}$. Thus $x_{2}^{a_{12}}x_{4}^{a_{34}}$ divides $-Q$, so $u_{2}-a_{2} \geq a_{12}$ a contradiction.

\begin{proposition} \label{Difficult2} Suppose that $I(C)$ is given as in case 2(a) and also that $a_{24}<a_{34}$. If $C$ has Cohen-Macaulay tangent cone at the origin then  \begin{enumerate} \item $a_{3} \leq a_{31}+a_{34}$, \item $a_{12}+a_{34} \leq a_{41} +a_{23}$ and \item $a_{2}+a_{12} \leq a_{1}+a_{23}-a_{13}+a_{24}$.
\end{enumerate}

\end{proposition}
\noindent \textbf{Proof.} From Proposition \ref{prop2(a)} we have that conditions (1) and (2) are true. Suppose first that $a_{13} \leq a_{23}$. Again from Proposition \ref{prop2(a)}, for $u_{2}=a_{2}+a_{12}$ and $u_{3}=u_{4}=0$, we deduce that there exists a monomial $M$ with $1 \in {\rm supp}(M)$ such that $x_{2}^{a_{2}+a_{12}}-M \in I(C)$ and $a_{2}+a_{12} \leq {\rm deg}(M)$. Since the binomial $x_{2}^{a_{2}+a_{12}}-M$ belongs to $I(C)$, we have that $x_{2}^{a_{2}+a_{12}}-M=\sum_{i=1}^{5}H_{i}f_{i}$ for some polynomials $H_{i} \in K[x_{1},\ldots,x_{4}]$. Then $x_{2}^{a_{2}+a_{12}}$ arise as a term in the sum $\sum_{i=1}^{5}H_{i}f_{i}$, so $Q=-x_{2}^{a_{12}}x_{3}^{a_{23}}x_{4}^{a_{24}}$ is a term in the sum $\sum_{i=1}^{5}H_{i}f_{i}$, therefore $Q$ should be canceled with another term of the above sum. Thus $T=-x_{1}^{a_{1}}x_{3}^{a_{23}-a_{13}}x_{4}^{a_{24}}$ is a term in the sum $\sum_{i=1}^{5}H_{i}f_{i}$. Suppose that $M \neq x_{1}^{a_{1}}x_{3}^{a_{23}-a_{13}}x_{4}^{a_{24}}$. Then $T$ should be canceled with another term of the above sum, a contradiction. Thus $M=x_{1}^{a_{1}}x_{3}^{a_{23}-a_{13}}x_{4}^{a_{24}}$ and also $a_{2}+a_{12} \leq a_{1}+a_{23}-a_{13}+a_{24}$.\\ Suppose now that $a_{13}>a_{23}$. Note that $x_{1}^{a_{1}}x_{4}^{a_{24}}-x_{2}^{a_{2}+a_{12}}x_{3}^{a_{13}-a_{23}}\in I(C)$. Thus there is a monomial $N$ with $1 \in {\rm supp}(N)$ such that $x_{2}^{a_{2}+a_{12}}x_{3}^{a_{13}-a_{23}}-N \in I(C)$ and $a_{2}+a_{12}+a_{13}-a_{23} \leq {\rm deg}(N)$. We have that $x_{2}^{a_{2}+a_{12}}x_{3}^{a_{13}-a_{23}}-N=\sum_{i=1}^{5}H_{i}f_{i}$ for some polynomials $H_{i} \in K[x_{1},\ldots,x_{4}]$. Then $x_{2}^{a_{2}+a_{12}}x_{3}^{a_{13}-a_{23}}$ arise as a term in the sum $\sum_{i=1}^{5}H_{i}f_{i}$, so $Q=-x_{2}^{a_{12}}x_{3}^{a_{13}}x_{4}^{a_{24}}$ is a term in the sum $\sum_{i=1}^{5}H_{i}f_{i}$, therefore $Q$ should be canceled with another term of the above sum. Thus $T=-x_{1}^{a_{1}}x_{4}^{a_{24}}$ is a term in the sum $\sum_{i=1}^{5}H_{i}f_{i}$. Let us suppose that $N \neq x_{1}^{a_{1}}x_{4}^{a_{24}}$. But then $T$ should be canceled with another term of the above sum a contradiction. Thus $N=x_{1}^{a_{1}}x_{4}^{a_{24}}$ and also $a_{2}+a_{12}+a_{13}-a_{23} \leq a_{1}+a_{24}$.

\begin{theorem} Suppose that $I(C)$ is given as in case 2(a) and also that $a_{24}<a_{34}$. Assume that $a_{13}\leq a_{23}$. Then $C$ has Cohen-Macaulay tangent cone at the origin if and only if \begin{enumerate} \item $a_{3} \leq a_{31}+a_{34}$, \item $a_{12}+a_{34} \leq a_{41} +a_{23}$ and \item $a_{2}+a_{12} \leq a_{1}+a_{23}-a_{13}+a_{24}$.
\end{enumerate}

\end{theorem}

\noindent \textbf{Proof.} $(\Longrightarrow)$ From Proposition \ref{Difficult2} we have that conditions (1), (2) and (3) are true.\\
$(\Longleftarrow)$ From Proposition \ref{prop2(a)} it is enough to consider a monomial $N=x_{2}^{u_2}x_{3}^{u_3}x_{4}^{u_4}$, where $u_2 \geq a_{2}$, $u_3<a_{13}$ and $u_4<a_{34}$, with the property: there exists at least one monomial $P$ such that $1 \in {\rm supp}(P)$ and also $N-P$ is in $I(C)$. Suppose that $u_{2} \geq a_{2}+a_{12}$. Let $P=x_{1}^{a_{1}}x_{2}^{u_{2}-a_{2}-a_{12}}x_{3}^{u_{3}+a_{23}-a_{13}}x_{4}^{u_{4}+a_{24}}$. We have that $N-P \in I(C)$ and $${\rm deg}(N)=u_{2}+u_{3}+u_{4} \leq u_{2}+u_{3}+u_{4}+a_{1}+a_{23}-a_{13}+a_{24}-a_{2}-a_{12}={\rm deg}(P).$$ It suffices to assume that $u_{2}-a_{2}<a_{12}$. Since the binomial $N-P$ belongs to $I(C)$, we have that $N-P=\sum_{i=1}^{5}H_{i}f_{i}$ for some polynomials $H_{i} \in K[x_{1},\ldots,x_{4}]$. Then $T=-x_{2}^{u_{2}-a_{2}}x_{3}^{u_{3}+a_{23}}x_{4}^{u_{4}+a_{24}}$ is a term in the sum $\sum_{i=1}^{5}H_{i}f_{i}$ and it should be canceled with another term of the above sum. Thus $x_{2}^{a_{12}}x_{4}^{a_{34}}$ divides $-T$, so $u_{2} -a_{2} \geq a_{12}$ a contradiction.\\

\begin{example} {\rm Consider $n_{1}=627$, $n_{2}=1546$, $n_{3}=1662$ and $n_{4}=3377$. The toric ideal $I(C)$ is minimally generated by the set $$G=\{x_1^{18}- x_{2}^{3}x_{3}^{4}, x_{2}^{25}- x_{3}^{7}x_{4}^{8}, x_3^{11}-x_{1}^{13}x_{4}^{3}, x_{4}^{11}-x_{1}^{5}x_{2}^{22}, x_{1}^{5}x_{3}^{7}-x_{2}^{3}x_{4}^{3}\}.$$ Here $a_{24}=8$, $a_{34}=3$, $a_{12}=3$, $a_{1}=18$, $a_{2}=25$ and $a_{13}=4<7=a_{23}$. Note that $a_{3}=11<14=a_{31}+a_{34}$ and $a_{12}+a_{34}=6<12=a_{41}+a_{23}$. We have that $a_{2}+a_{12}=28<29=a_{1}+a_{23}-a_{13}+a_{24}$. Thus $C$ has Cohen-Macaulay tangent cone at the origin. Remark that $x_{2}^{28}-x_{1}^{5}x_{3}^{14}x_{4}^{5} \in I(C)$, but ${\rm deg}(x_2^{28})=28>24={\rm deg}(x_{1}^{5}x_{3}^{14}x_{4}^{5})$.}

\end{example}

\begin{remark} {\rm Suppose that $I(C)$ is given as in case 2(b).\\ (1) It holds that $a_{1}>a_{12}+a_{13}$, $a_{4}<a_{41}+a_{43}$ and $a_{24}+a_{13}<a_{41} +a_{32}$.\\ (2) If $a_{24} \leq a_{34}$, then $x_{3}^{a_{3}+a_{13}}-x_{1}^{a_{41}}x_{2}^{2a_{32}}x_{4}^{a_{34}-a_{24}} \in I(C)$.\\(3) If $a_{12} \leq a_{32}$, then the binomial $x_{3}^{a_{3}+a_{13}}-x_{1}^{a_{1}}x_{2}^{a_{32}-a_{12}}x_{4}^{a_{34}}$ belongs to $I(C)$.}
\end{remark}

\begin{proposition} \label{Case 2(b) basic} Suppose that $I(C)$ is given as in case 2(b). Then $C$ has Cohen-Macaulay tangent cone at the origin if and only if \begin{enumerate} \item $a_{2} \leq a_{21}+a_{24}$ and \item for every monomial $x_{2}^{u_{2}}x_{3}^{u_{^3}}x_{4}^{u_{4}}$, where $u_{2}<a_{12}$, $u_{3} \geq a_{3}$ and $u_{4}<a_{24}$, with $u_{2}n_{2}+u_{3}n_{3}+u_{4}n_{4} \in n_{1}+\mathcal{S}$ there exists a monomial $N$ with $1 \in {\rm supp}(N)$ such that $M-N \in I(C)$ and also ${\rm deg}(M) \leq {\rm deg}(N)$.
\end{enumerate}

\end{proposition}

\noindent \textbf{Proof.} In this case $I(C)$ is minimally generated by the set $$G=\{f_{1} = x_1^{a_1}- x_2^{a_{12}} x_3^{a_{13}}, f_2 = x_{2}^{a_2}- x_{1}^{a_{21}}x_{4}^{a_{24}}, f_3 = x_3^{a_{3}}-x_{2}^{a_{32}}x_{4}^{a_{34}},$$ $$f_4 = x_{4}^{a_4}-x_{1}^{a_{41}}x_{3}^{a_{43}}, f_5 = x_{1}^{a_{41}}x_{2}^{a_{32}}-x_{3}^{a_{13}}x_4^{a_{24}}\}.$$
Suppose that $C$ has Cohen-Macaulay tangent cone at the origin. Since $I(C)$ is generated by the indispensable binomials, every binomial $f_i$, $1 \leq i \leq 5$, is indispensable of $I(C)$. In particular the binomial $f_2$ is indispensable of $I(C)$. Therefore the inequality $a_{2} \leq a_{21}+a_{24}$ holds. By Theorem \ref{BasicHer} condition (2) is also true.

Conversely, from Proposition \ref{1(b)} (iii), it is enough to consider a monomial $M=x_{2}^{u_2}x_{3}^{u_3}x_{4}^{u_4}$, where $u_2<a_{2}$, $u_3 \geq 0$ and $u_4<a_{4}$, with the property: there exists at least one monomial $P$ such that $1 \in {\rm supp}(P)$ and also ${\rm deg}_{\mathcal{S}}(M)={\rm deg}_{\mathcal{S}}(P)$. Suppose that $u_{3} \geq a_{3}$. If $u_{2} \geq a_{12}$, then we let $P=x_{1}^{a_{1}}x_{2}^{u_{2}-a_{12}}x_{3}^{u_{3}-a_{13}}x_{4}^{u_{4}}$. So we have that $M-P \in I(C)$ and also ${\rm deg}(M)<{\rm deg}(P)$. Similarly if $u_{4} \geq a_{24}$, then we let $P=x_{1}^{a_{41}}x_{2}^{u_{2}+a_{32}}x_{3}^{u_{3}-a_{13}}x_{4}^{u_{4}-a_{24}}$. So we have that $M-P \in I(C)$ and also ${\rm deg}(M)<{\rm deg}(P)$. If both conditions $u_{2}<a_{12}$ and $u_{4}<a_{24}$ hold, then condition (2) implies that there exists a monomial $N$ with $1 \in {\rm supp}(N)$ such that $M-N \in I(C)$ and also ${\rm deg}(M) \leq {\rm deg}(N)$. Suppose now that $u_3<a_3$. Then $M$ is divided by at least one of the monomials $x_{2}^{a_{12}}x_3^{a_{13}}$, $x_{3}^{a_{13}}x_4^{a_{24}}$ and $x_{2}^{a_{32}}x_4^{a_{34}}$. If $M$ is divided by $x_{2}^{a_{12}}x_3^{a_{13}}$, then $M=x_2^{a_{12}+p}x_3^{a_{13}+q}x_{4}^{u_4}$, for some non-negative integers $p$ and $q$, so $M-x_1^{a_{1}}x_2^{p}x_3^{q}x_{4}^{u_{4}} \in I(C)$ and also ${\rm deg}(M)<{\rm deg}(x_1^{a_{1}}x_2^{p}x_3^{q}x_{4}^{u_{4}})$. If $M$ is divided by $x_{3}^{a_{13}}x_4^{a_{24}}$, then $M=x_{2}^{u_2}x_3^{a_{13}+p}x_4^{a_{24}+q}$, for some non-negative integers $p$ and $q$, and therefore the binomial $M-x_1^{a_{41}}x_2^{u_{2}+a_{32}}x_3^{p}x_4^{q} \in I(C)$ and also ${\rm deg}(M) \leq {\rm deg}(x_1^{a_{41}}x_2^{u_{2}+a_{32}}x_3^{p}x_4^{q})$. Assume that neither $x_{2}^{a_{12}}x_3^{a_{13}}$ nor $x_{3}^{a_{13}}x_4^{a_{24}}$ divides $M$. Then necessarily $x_{2}^{a_{32}}x_4^{a_{34}}$ divides $M$. But $M$ is not divided by any leading monomial of $G$ with respect to $lex-inf$ such that $x_{1}>_{lex} x_{2} >_{lex} x_{3} >_{lex} x_{4}$. Thus $m=u_{2}n_{2}+u_{3}n_{3}+u_{4}n_{4}$ is in $Q$, a contradiction to the fact that $m-n_{1} \in \mathcal{S}$. Therefore, from Proposition \ref{1(b)}, $C$ has Cohen-Macaulay tangent cone at the origin.

The proof of the following proposition is similar to that of Proposition \ref{CasesCohen1(b)} and therefore it is omitted.

\begin{proposition} \label{CasesCohen2(b)} Suppose that $I(C)$ is given as in case 2(b). Assume that $C$ has Cohen-Macaulay tangent cone at the origin and also $a_{24} \leq a_{34}$. \begin{enumerate} \item If $a_{32}<a_{12}$, then $a_{3}+a_{13} \leq a_{41}+2a_{32}+a_{34}-a_{24}$.
\item If $a_{12} \leq a_{32}$, then $a_{3}+a_{13} \leq a_{1}+a_{32}-a_{12}+a_{34}$.
\end{enumerate}

\end{proposition}

\begin{theorem} Suppose that $I(C)$ is given as in case 2(b) and also that $a_{24} \leq a_{34}$. Then $C$ has Cohen-Macaulay tangent cone at the origin if and only if \begin{enumerate} \item $a_{2} \leq a_{21}+a_{24}$ and \item either $a_{32}<a_{12}$ and $a_{3}+a_{13} \leq a_{41}+2a_{32}+a_{34}-a_{24}$ or $a_{12} \leq a_{32}$ and $a_{3}+a_{13} \leq a_{1}+a_{32}-a_{12}+a_{34}$.
\end{enumerate}

\end{theorem}

\noindent \textbf{Proof.} $(\Longrightarrow)$ From Proposition \ref{Case 2(b) basic} we have that condition (1) is true. From Proposition \ref{CasesCohen2(b)} the condition (3) is also true.\\
$(\Longleftarrow)$ From Proposition \ref{Case 2(b) basic} it is enough to consider a monomial $N=x_{2}^{u_2}x_{3}^{u_3}x_{4}^{u_4}$, where $u_2<a_{12}$, $u_3 \geq a_{3}$ and $u_4<a_{24}$, with the property: there exists at least one monomial $P$ such that $1 \in {\rm supp}(P)$ and also $N-P$ is in $I(C)$. Suppose that $u_{3} \geq a_{3}+a_{13}$ and let $M$ denote either the monomial $x_{1}^{a_{41}}x_{2}^{2a_{32}}x_{4}^{a_{34}-a_{24}}$ when $a_{32}<a_{12}$ or the monomial $x_{1}^{a_{1}}x_{2}^{a_{32}-a_{12}}x_{4}^{a_{34}}$ when $a_{12} \leq a_{32}$. Let $P=x_{2}^{u_2}x_{3}^{u_{3}-a_{3}-a_{13}}x_{4}^{u_4}M$. We have that $N-P \in I(C)$ and $${\rm deg}(N)=u_{2}+u_{3}+u_{4} \leq u_{2}+u_{3}+u_{4}+{\rm deg}(M)-a_{3}-a_{13}={\rm deg}(P).$$ It suffices to consider the case that $u_{3}-a_{3}<a_{13}$. Since the binomial $N-P$ belongs to $I(C)$, we have that $N-P=\sum_{i=1}^{5}H_{i}f_{i}$ for some polynomials $H_{i} \in K[x_{1},\ldots,x_{4}]$. Now the monomial $N$ is divided by the monomial $x_{3}^{a_{3}}$, so $Q=-x_{2}^{u_{2}+a_{32}}x_{3}^{u_{3}-a_{3}}x_{4}^{u_{4}+a_{34}}$ is a term in the sum $\sum_{i=1}^{5}H_{i}f_{i}$ and it should be canceled with another term of the above sum. Remark that $u_{2}+a_{32}<a_{2}$ and $u_{4}+a_{34}<a_{4}$. Thus $x_{3}^{a_{13}}x_{4}^{a_{24}}$ divides $-Q$, so $u_{3}-a_{3} \geq a_{13}$ a contradiction.

\begin{proposition} Suppose that $I(C)$ is given as in case 2(b) and also that $a_{34}<a_{24}$. If $C$ has Cohen-Macaulay tangent cone then \begin{enumerate} \item $a_{2} \leq a_{21}+a_{24}$ and \item $a_{3}+a_{13} \leq a_{1}+a_{32}-a_{12}+a_{34}$.
\end{enumerate}

\end{proposition}

\begin{theorem} Suppose that $I(C)$ is given as in case 2(b) and also that $a_{34}<a_{24}$. Assume that $a_{12} \leq a_{32}$. Then $C$ has Cohen-Macaulay tangent cone at the origin if and only if \begin{enumerate} \item $a_{2} \leq a_{21}+a_{24}$ and \item $a_{3}+a_{13} \leq a_{1}+a_{32}-a_{12}+a_{34}$.
\end{enumerate}

\end{theorem}

\begin{example} {\rm Consider $n_{1}=813$, $n_{2}=1032$, $n_{3}=1240$ and $n_{4}=1835$. The toric ideal $I(C)$ is minimally generated by the set $$G=\{x_1^{16}- x_2^{9}x_{3}^{3}, x_{2}^{14}- x_{1}^{11}x_{4}^{3}, x_3^{16}-x_{2}^{5}x_{4}^{8}, x_{4}^{11}-x_{1}^{5}x_{3}^{13}, x_{1}^{5}x_{2}^{5}-x_{3}^{3}x_{4}^{3}\}.$$ Here $a_{13}=a_{24}=3$, $a_{34}=8$, $a_{41}=5$, $a_{3}=16$ and $a_{32}=5<9=a_{12}$. Note that $a_{2}=14=a_{21}+a_{24}$. We have that $a_{3}+a_{13}=19<20=a_{41}+2a_{32}+a_{34}-a_{24}$. Consequently $C$ has Cohen-Macaulay tangent cone at the origin.}

\end{example}

\begin{theorem} Suppose that $I(C)$ is given as in case 3(a). Then $C$ has Cohen-Macaulay tangent cone at the origin if and only if $a_{2} \leq a_{21}+a_{23}$ and $a_{3} \leq a_{31}+a_{34}$.

\end{theorem}

\textbf{Proof.} In this case $I(C)$ is minimally generated by the set $$G=\{f_{1} = x_1^{a_1}- x_2^{a_{12}} x_4^{a_{14}}, f_2 = x_{2}^{a_2}- x_{1}^{a_{21}}x_{3}^{a_{23}}, f_3 = x_3^{a_{3}}-x_{1}^{a_{31}}x_{4}^{a_{34}},$$ $$f_4 = x_{4}^{a_4}-x_{2}^{a_{42}}x_{3}^{a_{43}}, f_5 = x_{1}^{a_{31}}
x_{2}^{a_{42}}-x_{3}^{a_{23}}x_4^{a_{14}}\}.$$ If $a_{2} \leq a_{21}+a_{23}$ and $a_{3} \leq a_{31}+a_{34}$, then we have, from Theorem 2.10 in \cite{ArMe}, that the curve $C$ has Cohen-Macaulay tangent cone at the origin. Conversely suppose that $C$ has Cohen-Macaulay tangent cone at the origin. Since $I(C)$ is generated by the indispensable binomials, every binomial $f_i$, $1 \leq i \leq 5$, is indispensable of $I(C)$. In particular the binomials $f_2$ and $f_3$ are indispensable of $I(C)$. If there exists a monomial $N \neq x_{1}^{a_{31}}x_{4}^{a_{34}}$ such that $g=x_{3}^{a_3}-N$ belongs to $I(C)$, then we can replace $f_{3}$ in $S$ by the binomials $g$ and $N-x_{1}^{a_{31}}x_{4}^{a_{34}} \in I(C)$, a contradiction to the fact that $f_3$ is indispensable. Thus $N=x_{1}^{a_{31}}x_{4}^{a_{34}}$ and therefore, from Theorem \ref{BasicHer}, we have that $a_{3} \leq a_{31}+a_{34}$. Similarly we get that $a_{2} \leq a_{21}+a_{23}$.

\begin{remark} {\rm Suppose that $I(C)$ is given as in case 3(b).\\ (1) It holds that $a_{1}>a_{12}+a_{14}$, $a_{2}>a_{23}+a_{24}$, $a_{3}<a_{31}+a_{32}$ and $a_{4}<a_{41}+a_{43}$.\\ (2) If $a_{43} \leq a_{23}$, then $x_{2}^{a_{2}+a_{12}}-x_{1}^{a_{31}}x_{3}^{a_{23}-a_{43}}x_{4}^{2a_{24}} \in I(C)$. \\(3) If $a_{14} \leq a_{24}$, then the binomial $x_{2}^{a_{2}+a_{12}}-x_{1}^{a_{1}}x_{3}^{a_{23}}x_{4}^{a_{24}-a_{14}}$ belongs to $I(C)$.}
\end{remark}

\begin{proposition} \label{prop3(b)} Suppose that $I(C)$ is given as in case 3(b). Then $C$ has Cohen-Macaulay tangent cone at the origin if and only if \begin{enumerate} \item $a_{12}+a_{43} \leq a_{31} +a_{24}$ and
\item for every monomial $M=x_{2}^{u_{2}}x_{3}^{u_{^3}}x_{4}^{u_{4}}$, where $u_{2} \geq a_{2}$, $u_{3}<a_{43}$ and $u_{4}<a_{14}$, with $u_{2}n_{2}+u_{3}n_{3}+u_{4}n_{4} \in n_{1}+\mathcal{S}$ there exists a monomial $N$ with $1 \in {\rm supp}(N)$ such that $M-N \in I(C)$ and also ${\rm deg}(M) \leq {\rm deg}(N)$.
\end{enumerate}

\end{proposition}

\textbf{Proof.} In this case $I(C)$ is minimally generated by the set $$G=\{f_{1} = x_1^{a_1}- x_2^{a_{12}} x_4^{a_{14}}, f_2 = x_{2}^{a_2}- x_{3}^{a_{23}}x_{4}^{a_{24}}, f_3 = x_3^{a_{3}}-x_{1}^{a_{31}}x_{2}^{a_{32}},$$ $$f_4 = x_{4}^{a_4}-x_{1}^{a_{41}}x_{3}^{a_{43}}, f_5 = x_{1}^{a_{31}}x_{4}^{a_{24}}-x_{2}^{a_{12}}x_3^{a_{43}}\}.$$
Suppose that $C$ has Cohen-Macaulay tangent cone at the origin. Since $I(C)$ is generated by the indispensable binomials, every binomial $f_i$, $1 \leq i \leq 5$, is indispensable of $I(C)$. In particular the binomial $f_5$ is indispensable of $I(C)$. Therefore the inequality $a_{12}+a_{43} \leq a_{31} +a_{24}$ holds. By Theorem \ref{BasicHer} condition (2) is also true.

Conversely, from Proposition \ref{2(a)} (i), it is enough to consider a monomial $M=x_{2}^{u_2}x_{3}^{u_3}x_{4}^{u_4}$, where $u_2 \geq a_{2}$, $u_3<a_{3}$ and $u_4<a_{4}$, with the property: there exists at least one monomial $P$ such that $1 \in {\rm supp}(P)$ and also $M-P$ is in $I(C)$. If $u_{3} \geq a_{43}$, then we let $P=x_{1}^{a_{31}}x_{2}^{u_{2}-a_{12}}x_{3}^{u_{3}-a_{43}}x_{4}^{u_{4}+a_{24}}$. So we have that $M-P \in I(C)$ and also ${\rm deg}(M) \leq {\rm deg}(P)$. Similarly if $u_{4} \geq a_{14}$, then we let $P=x_{1}^{a_{1}}x_{2}^{u_{2}-a_{12}}x_{3}^{u_{3}}x_{4}^{u_{4}-a_{14}}$. So we have that $M-P \in I(C)$ and also ${\rm deg}(M)<{\rm deg}(P)$. If both conditions $u_{3}<a_{43}$ and $u_{4}<a_{14}$ hold, then condition (2) implies that there exists a monomial $N$ with $1 \in {\rm supp}(N)$ such that $M-N \in I(C)$ and also ${\rm deg}(M) \leq {\rm deg}(N)$. Therefore, from Proposition \ref{2(a)}, $C$ has Cohen-Macaulay tangent cone at the origin.

The proof of the next proposition is similar to that of Proposition \ref{CasesCohen1(b)} and therefore it is omitted.

\begin{proposition} \label{CasesCohen3(b)} Suppose that $I(C)$ is given as in case 3(b). Assume that $C$ has Cohen-Macaulay tangent cone at the origin and also $a_{43} \leq a_{23}$. \begin{enumerate} \item If $a_{24}<a_{14}$, then $a_{2}+a_{12} \leq a_{31}+2a_{24}+a_{23}-a_{43}$.
\item If $a_{14} \leq a_{24}$, then $a_{2}+a_{12} \leq a_{1}+a_{23}+a_{24}-a_{14}$.
\end{enumerate}

\end{proposition}

\begin{theorem} Suppose that $I(C)$ is given as in case 3(b) and $a_{43} \leq a_{23}$. Then $C$ has Cohen-Macaulay tangent cone at the origin if and only if \begin{enumerate} \item $a_{12}+a_{43} \leq a_{31} +a_{24}$
and \item either $a_{24}<a_{14}$ and $a_{2}+a_{12} \leq a_{31}+2a_{24}+a_{23}-a_{43}$ or $a_{14} \leq a_{24}$ and $a_{2}+a_{12} \leq a_{1}+a_{23}+a_{24}-a_{14}$.
\end{enumerate}

\end{theorem}

\noindent \textbf{Proof.} $(\Longrightarrow)$ From Proposition \ref{prop3(b)} we have that conditions (1) and (2) are true. From Proposition \ref{CasesCohen3(b)} the condition (3) is also true.\\
$(\Longleftarrow)$ From Proposition \ref{prop3(b)} it is enough to consider a monomial $N=x_{2}^{u_2}x_{3}^{u_3}x_{4}^{u_4}$, where $u_2 \geq a_{2}$, $u_3<a_{43}$ and $u_4<a_{14}$, with the property: there exists at least one monomial $P$ such that $1 \in {\rm supp}(P)$ and also $N-P$ is in $I(C)$. Suppose that $u_{2} \geq a_{2}+a_{12}$ and let $M$ denote either the monomial $x_{1}^{a_{31}}x_{3}^{a_{23}-a_{43}}x_{4}^{2a_{24}}$ when $a_{24}<a_{14}$ or the monomial $x_{1}^{a_{31}}x_{3}^{a_{23}}x_{4}^{a_{24}-a_{14}}$ when $a_{14} \leq a_{24}$. Let $P=x_{2}^{u_{2}-a_{2}-a_{12}}x_{3}^{u_{3}}x_{4}^{u_4}M$. We have that $N-P \in I(C)$ and $${\rm deg}(N)=u_{2}+u_{3}+u_{4} \leq u_{2}+u_{3}+u_{4}+{\rm deg}(M)-a_{2}-a_{12}={\rm deg}(P).$$ It suffices to consider the case that $u_{2}-a_{2}<a_{12}$. Since the binomial $N-P$ belongs to $I(C)$, we have that $N-P=\sum_{i=1}^{5}H_{i}f_{i}$ for some polynomials $H_{i} \in K[x_{1},\ldots,x_{4}]$. Now the monomial $N$ is divided by the monomial $x_{2}^{a_{2}}$, so $Q=-x_{2}^{u_{2}-a_{2}}x_{3}^{u_{3}+a_{23}}x_{4}^{u_{4}+a_{24}}$ is a term in the sum $\sum_{i=1}^{5}H_{i}f_{i}$ and it should be canceled with another term of the above sum. Remark that $u_{3}+a_{23}<a_{3}$ and $u_{4}+a_{24}<a_{4}$. Thus $x_{2}^{a_{12}}x_{3}^{a_{43}}$ divides $-Q$, so $u_{2}-a_{2} \geq a_{12}$ a contradiction.

\begin{proposition} Suppose that $I(C)$ is given as in case 3(b) and also that $a_{23}<a_{43}$. If $C$ has Cohen-Macaulay tangent cone at the origin then \begin{enumerate} \item $a_{12}+a_{43} \leq a_{31} +a_{24}$ and \item $a_{2}+a_{12} \leq a_{1}+a_{23}+a_{24}-a_{14}$.
\end{enumerate}

\end{proposition}

\begin{theorem} Suppose that $I(C)$ is given as in case 3(b) and also that $a_{23}<a_{43}$. Assume that $a_{14} \leq a_{24}$. Then $C$ has Cohen-Macaulay tangent cone at the origin if and only if \begin{enumerate} \item $a_{12}+a_{43} \leq a_{31} +a_{24}$ and \item $a_{2}+a_{12} \leq a_{1}+a_{23}+a_{24}-a_{14}$.
\end{enumerate}

\end{theorem}

\section{Families of monomial curves supporting Rossi's problem}

In this section, we give some examples showing how the criteria given in the previous one can be used to give families of monomial curves supporting Rossi's problem in $\mathbb{A}^{4}(K)$.

\begin{example}
{\rm Consider the family  	
$n_{1}=m^{3}+m^{2}-m$, $n_{2}=m^{3}+2m^{2}+m-1$, $n_{3}=m^{3}+3m^{2}+2m-2$ and $n_{4}=m^{3}+4m^{2}+3m-2$
for $m \geq 2$ given in \cite{ArMe}. Let $C_{m}$ be the corresponding monomial curve.
The toric ideal $I(C_{m})$ is generated by the set
$$S_{m}=\{x_{1}^{m+3}-x_{3}x_{4}^{m-1},
x_{2}^{m+2}-x_{1}^{m+2}x_{4},
x_{3}^{m}-x_{1}x_{2}^{m},
x_{4}^{m}-x_{2}^{2}x_{3}^{m-1},$$ $$x_{1}^{m+2}x_{3}^{m-1}-x_{2}^{m}x_{4}^{m-1}\}.$$
Thus we are in case 1(a) of Remark \ref{BasicGorenstein} and it is sufficient to consider the binomial $x_{2}^{m+2}-x_{1}^{m+2}x_{4}$, which guarantees that $C_{m}$ has Cohen-Macaulay
tangent cone. For each fixed $m$, by using the technique given in \cite{Gi}, we can construct a new family of monomial curves having Cohen-Macaulay tangent cone. For $m=2$, we have the symmetric semigroup generated by $n_{1}=10$, $n_{2}=17$, $n_{3}=22$ and $n_{4}=28$. The corresponding monomial curve is $C_{2}$ and $I(C_{2})$ is minimally generated by
the set
$$S_2=\{ x_{1}^{5}-x_{3}x_{4}, x_{2}^{4}-x_{1}^{4}x_{4}, x_{3}^{2}-x_{1}x_{2}^{2}, x_{4}^{2}-x_{2}^{2}x_{3}, x_{1}^{4}x_{3}-x_{2}^{2}x_{4}\}.$$
By the method given in \cite{Gi}, the semigroup generated by  $m_{1}=10+6t$, $m_{2}=17+9t$, $m_{3}=22+6t$ and $m_{4}=28+12t$ (for $t$ a non-negative integer) is symmetric, whenever ${\rm gcd}(10+6t,17+9t,22+6t, 28+12t)=1$. Moreover, the toric ideal of the corresponding monomial curve is minimally generated by the set
$$\{ x_{1}^{t+5}-x_{3}^{t+1}x_{4}, x_{2}^{4}-x_{1}^{4}x_{4}, x_{3}^{t+2}-x_{1}^{t+1}x_{2}^{2}, x_{4}^{2}-x_{2}^{2}x_{3}, x_{1}^{4}x_{3}-x_{2}^{2}x_{4}\}$$
by the construction in \cite{Gi}.
Here, if $t=1$, then $m_1< m_2 < m_3 < m_4$, and we are in case 1(a) again. From the binomial $x_{2}^{4}-x_{1}^{4}x_{4}$, we deduce that the corresponding monomial curve has Cohen-Macaulay tangent cone. If $t \geq 1$, then $m_2 > m_3$. In this case, we interchange them to get the semigroup generated by  $m'_{1}=10+6t$, $m'_{2}=22+6t$, $m'_{3}=17+9t$ and $m'_{4}=28+12t$ with $m'_1< m_2' < m_3' < m_4'$. Thus the toric ideal of the corresponding monomial curve is generated by the set
$$\{ x_{1}^{t+5}-x_{2}^{t+1}x_{4}, x_{2}^{t+2}-x_{1}^{t+1}x_{3}^{2}, x_{3}^{4}-x_{1}^{4}x_{4}, x_{4}^{2}-x_{2}x_{3}^{2}, x_{1}^{4}x_{2}-x_{3}^{2}x_{4}\}.$$
Now we are in case 3(a) of Remark \ref{BasicGorenstein} and the binomials $x_{3}^{4}-x_{1}^{4}x_{4}$ and $x_{2}^{t+2}-x_{1}^{t+1}x_{3}^{2}$ guarantee that the corresponding monomial curve has Cohen-Macaulay tangent cone. In this way, we can construct infinitely many families of Gorenstein monomial curves having Cohen-Macaulay tangent cones. In other words, the corresponding local rings have non-decreasing Hilbert functions supporting Rossi's problem.}

\end{example}

In the literature, there are no examples of non-complete intersection Gorenstein monomial curve families supporting Rossi's problem, although their tangent cones are not Cohen-Macaulay. The next example gives a family of monomial curves with the above property. To prove it we will use the following proposition. 
\begin{proposition} \cite[Proposition 2.2]{bayer} \label{hilbert} Let $I \subset K[x_{1}, x_{2}, \ldots, x_{d}]$ be a monomial ideal and $I=\langle J, {\bf x}^{\bf u} \rangle$ for a monomial ideal $J$ and a monomial ${\bf x}^{\bf u}$. Let $p(I)$ denote the numerator $g(t)$ of the Hilbert Series for $K[x_{1}, x_{2}, \ldots, x_{d}]/I$, and let $|{\bf u}|$ denote the total degree of the monomial ${\bf x}^{\bf u}$. Then $p(I)=p(J)-t^{|{\bf u}|}p(J:{\bf x}^{\bf u})$.
\end{proposition}
	
\begin{example}

{\rm Consider the family
$n_{1}=2m+1$, $n_{2}=2m+3$, $n_{3}=2m^{2}+m-2$ and $n_{4}=2m^2+m-1$, where $m \geq 4$ is an integer. Let $C_{m}$ be the corresponding monomial curve. The toric ideal $I(C_{m})$ is minimally generated by the binomials
$$x_{1}^{m+1}-x_{2}x_{3}, x_{2}^{m}-x_{1}x_{4}, x_{3}^{2}-x_{2}^{m-1}x_{4}, x_{4}^{2}-x_{1}^{m}x_{3}, x_{1}^{m}x_{2}^{m-1}-x_{3}x_{4}.$$
Thus we are in Case 2(b) of Remark \ref{BasicGorenstein}. Consider the binomial $x_{2}^{m}-x_{1}x_{4}$. Since $m \geq 4$, we have, from Theorem \ref{Case 2(b) basic}, that the tangent cone of $C_{m}$ is not Cohen-Macaulay. It is enough to show that the Hilbert function of $K[x_{1}, x_{2}, x_3, x_{4}]/I(C_{m})_*$ is non-decreasing. Where $I(C_{m})_*$ is the ideal generated by the polynomials $f_*$ for $f$ in $I(C_{m})$ and $f_*$ is the homogeneous summand of $f$ of least degree. By a standard basis computation, $I(C_{m})_*$ is generated by the set
$$\{x_{2}x_{3}, x_{3}^{2}, x_{1}x_{4}, x_{3}x_{4}, x_{4}^{2}, x_{2}^{m}x_{4}, x_{1}^{m+2}x_{3}, x_{2}^{2m+1}\}.$$
Let $$J_{0}=I(C_m)_*, J_{1}=\langle x_{3}^{2}, x_{1}x_{4}, x_{3}x_{4}, x_{4}^{2}, x_{2}^{m}x_{4}, x_{1}^{m+2}x_{3}, x_{2}^{2m+1} \rangle,$$ $$J_{2}=\langle x_{3}^{2}, x_{3}x_{4}, x_{4}^{2}, x_{2}^{m}x_{4}, x_{1}^{m+2}x_{3}, x_{2}^{2m+1} \rangle, J_{3}=\langle x_{3}^{2}, x_{4}^{2}, x_{2}^{m}x_{4}, x_{1}^{m+2}x_{3}, x_{2}^{2m+1} \rangle.$$ Note that $J_i= \langle J_{i+1}, q_i \rangle$, where $q_0=x_2x_3$, $q_1=x_1x_4$ and $q_2=x_3x_4$. We apply Proposition \ref{hilbert} to the ideal $J_i$ for $0 \leq i \leq 2$, so
\begin{equation} \label{recursion}
p(J_i)=p(J_{i+1})-t^{2}p(J_{i+1}:q_i).
\end{equation}
In this case, we have
$J_{1}:(x_2x_3)=\langle x_3,x_4,x_{1}^{m+2}, x_2^{2m} \rangle$,
$J_{2}: (x_1x_4)=\langle x_3,x_4, x_2^{m} \rangle$ and $J_{3}: (x_3x_4)= \langle x_3, x_4, x_{1}^{m+2}, x_2^m \rangle$. Since
$$K[x_{1}, x_{2}, x_3, x_{4}]/ \langle x_{3}^{2}, x_{1}^{m+2}x_{3},
x_{4}^{2}, x_{2}^{m}x_{4},x_{2}^{2m+1} \rangle$$ is isomorphic to
$$K[x_1, x_3]/ \langle x_{3}^{2},x_{1}^{m+2}x_{3} \rangle \otimes K[x_{2},x_{4}]/ \langle x_{4}^{2}, x_{2}^{m}x_{4},x_{2}^{2m+1} \rangle,$$
we obtain
$$p(J_3)=(1-t)^{3}(1+t-t^{m+3})(1+2t+\cdots+2t^m+t^{m+1}+\cdots
+t^{2m}).$$ Substituting all these recursively in Equation (\ref{recursion}),
we obtain that the Hilbert series of
$K[x_{1}, x_{2}, x_3, x_{4}]/J_0$ is
$$\frac{1+3t+t^2+t^3+\cdots+t^m+t^{m+2}+t^{m+4}+t^{m+5}+
\cdots+t^{2m}}{1-t}.$$
Since the numerator does not have any negative coefficients, the Hilbert function is non-decreasing. In this way, we have shown that the Hilbert function of the local ring corresponding to the non-complete intersection Gorenstein monomial curve $C_m$ is non-decreasing for $m \geq 4$.}
\end{example}

\section{Appendix}

In the appendix we provide results concerning all the cases of Theorem \ref{KaOjVeryBasic}, except from Case 1. The proofs are similar to those in section 2 and therefore omitted. Let $n_{1}<n_{2}<n_{3}<n_{4}$ be positive integers with ${\rm gcd}(n_{1},\ldots,n_{4})=1$.

\begin{theorem} Suppose that $I(C)$ is given as in case 2(a). Let $\{x_1^{a_1} -
x_i^{a_i}, x_j^{a_j} - x_k^{a_k}, x_1^{a_1} -
x_j^{u_j}x_k^{u_k}\}$ be a generating set of $\mathcal{C}_{A}$ and assume that $n_{k}>n_{j}$. Then $C$ has Cohen-Macaulay tangent cone at the origin if and only if \begin{enumerate} \item for every binomial $f=M-N \in \mathcal{I}$ with $1 \in {\rm supp}(M)$ we have that ${\rm deg}(N) \leq {\rm deg}(M)$, \item for every binomial $f=M-N \in \mathcal{R}$ with $1 \in {\rm supp}(M)$ there is a monomial $P$ with $1 \in {\rm supp}(P)$ such that ${\rm deg}_{\mathcal{S}}(N)={\rm deg}_{\mathcal{S}}(P)$ and ${\rm deg}(N) \leq {\rm deg}(P)$ and \item for every monomial $M=x_{i}^{v_{i}}x_{j}^{v_{^j}}x_{k}^{v_{k}}$, where $v_{i}<a_{i}$, $v_{j} \geq a_{j}$ and $v_{k}<a_{k}$, with $v_{i}n_{i}+v_{j}n_{j}+v_{k}n_{k} \in n_{1}+\mathcal{S}$ there exists a monomial $N$ with $1 \in {\rm supp}(N)$ such that $M-N \in I(C)$ and also ${\rm deg}(M) \leq {\rm deg}(N)$.
\end{enumerate}

\end{theorem}

\begin{theorem} \label{Th2(a)} Suppose that $I(C)$ is given as in case 2(a). Let $\{x_1^{a_1} -x_i^{a_i}, x_j^{a_j} - x_k^{a_k}, x_j^{a_j}-x_1^{u_1}x_i^{u_i}\}$ be a generating set of $\mathcal{C}_{A}$ and assume that $n_{k}>n_{j}$. Then $C$ has Cohen-Macaulay tangent cone at the origin if and only if \begin{enumerate} \item for every binomial $f=M-N \in \mathcal{I}$ with $1 \in {\rm supp}(M)$ we have that ${\rm deg}(N) \leq {\rm deg}(M)$, \item for every $f=M-N \in \mathcal{R}$ with $1 \in {\rm supp}(M)$ there is a monomial $P$ with $1 \in {\rm supp}(P)$ such that ${\rm deg}_{\mathcal{S}}(N)={\rm deg}_{\mathcal{S}}(P)$ and ${\rm deg}(N) \leq {\rm deg}(P)$ and \item for every monomial $M=x_{i}^{v_{i}}x_{j}^{v_{j}}x_{k}^{v_{k}}$, where $v_{i}<a_{i}$, $v_{j} \geq a_{j}$ and $v_{k}<a_{k}$, there exists a monomial $N$ with $1 \in {\rm supp}(N)$ such that $M-N \in I(C)$ and also ${\rm deg}(M) \leq {\rm deg}(N)$.
\end{enumerate}

\end{theorem}

\begin{corollary} \label{BasicCorollary} Suppose that $I(C)$ is given as in case 2(a). Let $\{x_1^{a_1} -
x_i^{a_i}, x_j^{a_j} - x_k^{a_k}, x_j^{a_j} -
x_1^{u_1}x_i^{u_i}\}$ be a generating set of $\mathcal{C}_{A}$ and assume that $n_{k}>n_{j}$. If \begin{enumerate} \item for every binomial $f=M-N \in \mathcal{I}$ with $1 \in {\rm supp}(M)$ we have that ${\rm deg}(N) \leq {\rm deg}(M)$, \item for every $f=M-N \in \mathcal{R}$ with $1 \in {\rm supp}(M)$ there is a monomial $P$ with $1 \in {\rm supp}(P)$ such that ${\rm deg}_{\mathcal{S}}(N)={\rm deg}_{\mathcal{S}}(P)$ and ${\rm deg}(N) \leq {\rm deg}(P)$ and \item $a_{j} \leq u_{1}+u_{i}$,
\end{enumerate}
then $C$ has Cohen-Macaulay tangent cone at the origin.
\end{corollary}

\begin{example} {\rm Consider $n_{1}=30$, $n_{2}=34$, $n_{3}=42$ and $n_{4}=51$. The toric ideal $I(C)$ is minimally generated by the binomials $f_{1} = x_1^{7}- x_3^{5}$, $f_2 = x_{2}^{3}-x_{4}^{2}$ and $f_3 = x_2^{3}-x_{1}^{2}x_{3}$. Note that there is no binomial $g \in \mathcal{I}$ and also no binomial $h \in \mathcal{R}$. Then $C$ has Cohen-Macaulay tangent cone by Corollary \ref{BasicCorollary}.}

\end{example}

\begin{theorem} \label{The case 2(b)} Suppose that $I(C)$ is given as in case 2(b). Let $\{x_1^{a_1} -
x_i^{a_i}, x_j^{a_j} - x_k^{a_k}\}$ be a generating set of $\mathcal{C}_{A}$ and assume that $n_{k}>n_{j}$. Then $C$ has Cohen-Macaulay tangent cone at the origin if and only if \begin{enumerate} \item for every binomial $f=M-N \in \mathcal{I}$ with $1 \in {\rm supp}(M)$ we have that ${\rm deg}(N) \leq {\rm deg}(M)$, \item for every $f=M-N \in \mathcal{R}$ with $1 \in {\rm supp}(M)$ there is a monomial $P$ with $1 \in {\rm supp}(P)$ such that ${\rm deg}_{\mathcal{S}}(N)={\rm deg}_{\mathcal{S}}(P)$ and ${\rm deg}(N) \leq {\rm deg}(P)$ and \item for every monomial $M=x_{i}^{u_{i}}x_{j}^{u_{j}}x_{k}^{u_{k}}$, where $u_{i}<a_{i}$, $u_{j} \geq a_{j}$ and $u_{k}<a_{k}$, with $u_{i}n_{i}+u_{j}n_{j}+u_{k}n_{k} \in n_{1}+\mathcal{S}$ there exists a monomial $N$ with $1 \in {\rm supp}(N)$ such that $M-N \in I(C)$ and also ${\rm deg}(M) \leq {\rm deg}(N)$.
\end{enumerate}

\end{theorem}

\begin{example} {\rm Consider $n_{1}=49$, $n_{2}=63$, $n_{3}=65$ and $n_{4}=78$. The toric ideal $I(C)$ is minimally generated by the binomials $f_{1} = x_1^{9}- x_2^{7}$, $f_2 = x_{3}^{6}-x_{4}^{5}$, $f_{3}=x_{1}^{2}x_{4}^{2}-x_{2}^{3}x_{3}$, $f_{4}=x_{1}^{3}x_{2}^{2}-x_{3}^{3}x_{4}$, $f_{5}=x_{1}^{5}x_{4}-x_{2}x_{3}^{4}$, $f_{6}=x_{1}x_{2}^{5}-x_{3}^{2}x_{4}^{3}$, $f_{7}=x_{1}^{5}x_{3}^{2}-x_{2}x_{4}^{4}$, $f_{8}=x_{1}^{8}x_{2}-x_{3}x_{4}^{5}$, $f_{9}=x_{1}^{2}x_{3}^{5}-x_{2}^{3}x_{4}^{3}$, $f_{10}=x_{1}^{7}x_{3}-x_{2}^{4}x_{4}^{2}$ and $f_{11}=x_{1}^{4}x_{3}^{4}-x_{2}^{6}x_{4}$. Here $a_{1}=9$, $a_{2}=7$, $a_{3}=6$ and $a_{4}=5$. Also $\mathcal{I}=\{f_{3},f_{4},f_{5},f_{6},f_{7},f_{9},f_{10},f_{11}\}$ and $\mathcal{R}=\{f_{8}\}$. Let $M=x_{2}^{u_{2}}x_{3}^{u_{3}}x_{4}^{u_{4}}$, where $0 \leq u_{2}<7$, $u_{3} \geq 6$ and $0 \leq u_{4}<5$, with $u_{2}n_{2}+u_{3}n_{3}+u_{4}n_{4} \in n_{1}+\mathcal{S}$. Note that $x_{3}^{7}-x_{1}^{8}x_{2} \in I(C)$, so for $u_{3} \geq 7$ we have that $x_{2}^{u_{2}}x_{3}^{u_{3}}x_{4}^{u_{4}}-x_{1}^{8}x_{2}^{u_{2}+1}x_{3}^{u_{3}-7}x_{4}^{u_{4}} \in I(C)$ and also $u_{2}+u_{3}+u_{4}<u_{2}+u_{3}+u_{4}+2$. Furthermore $x_{2}x_{3}^{6}-x_{1}^{5}x_{3}^{2}x_{4} \in I(C)$ and $x_{3}^{6}x_{4}-x_{1}^{3}x_{2}^{2}x_{3}^{3} \in I(C)$. Thus for $0 \leq u_{2}<7$, $u_{3} \geq 6$, $0 \leq u_{4}<5$ there exists a monomial $N$ with $1 \in {\rm supp}(N)$ such that $M-N \in I(C)$ and also ${\rm deg}(M) \leq {\rm deg}(N)$. From Theorem \ref{The case 2(b)} the monomial curve $C$ has Cohen-Macaulay tangent cone at the origin.}

\end{example}

\begin{remark} {\rm In the case 2(c) we have that $I(C)$ is generated by the set $\{x_1^{a_1} - x_i^{a_i}, x_i^{a_i} - x_j^{a_j}, x_j^{a_j}-x_k^{a_k}\}$, so $I(C)$ is complete intersection and also $C$ has Cohen-Macaulay tangent cone at the origin.}
\end{remark}

\begin{theorem} Suppose that $I(C)$ is given as in case 3. Let $\{x_i^{a_i} -
x_j^{a_j}, x_j^{a_j} - x_k^{a_k}, x_1^{a_1} -
x_i^{u_{i}}x_j^{u_{j}}\}$ be a generating set of $\mathcal{C}_{A}$ and $n_{i}<min\{n_{j},n_{k}\}$.\\ Then $C$ has Cohen-Macaulay tangent cone at the origin if and only if \begin{enumerate} \item for every binomial $f=M-N \in \mathcal{I}$ with $1 \in {\rm supp}(M)$ we have that ${\rm deg}(N) \leq {\rm deg}(M)$ and \item for every monomial $M=x_{i}^{v_{i}}x_{j}^{v_{j}}x_{k}^{v_{k}}$, where $v_{i} \geq a_{i}$, with $v_{i}n_{i}+v_{j}n_{j}+v_{k}n_{k} \in n_{1}+\mathcal{S}$ there exists a monomial $N$ with $1 \in {\rm supp}(N)$ such that $M-N \in I(C)$ and also ${\rm deg}(M) \leq {\rm deg}(N)$.
\end{enumerate}

\end{theorem}

\begin{remark} {\rm Suppose that $I(C)$ is given as in case 3. Assume that a generating set of $\mathcal{C}_{A}$ is either $\{x_2^{a_2} -
x_i^{a_i}, x_i^{a_i} - x_j^{a_j}, x_1^{a_1} -
x_i^{u_{i}}x_j^{u_{j}}\}$ or $\{x_2^{a_2} -
x_i^{a_i}, x_i^{a_i} - x_j^{a_j}, x_1^{a_1} -
x_{2}^{u_{2}}x_i^{u_{i}}x_j^{u_{j}}\}$ where $u_{2}$, $u_{i}$ and $u_{j}$ are positive integers. Then $C$ has Cohen-Macaulay tangent cone at the origin if and only if \begin{enumerate} \item for every binomial $f=M-N \in \mathcal{I}$ with $1 \in {\rm supp}(M)$ we have that ${\rm deg}(N) \leq {\rm deg}(M)$ and \item for every monomial $x_{2}^{v_{2}}x_{3}^{v_{3}}x_{4}^{v_{4}}$, where $v_{2} \geq a_{2}$, with $v_{2}n_{2}+v_{3}n_{3}+v_{4}n_{4} \in n_{1}+\mathcal{S}$ there exists a monomial $N$ with $1 \in {\rm supp}(N)$ such that $M-N \in I(C)$ and also ${\rm deg}(M) \leq {\rm deg}(N)$.
\end{enumerate}}

\end{remark}
\begin{theorem} Suppose that $I(C)$ is given as in case 3. Let $\{x_1^{a_1} -
x_i^{a_i}, x_i^{a_i} - x_j^{a_j}, x_k^{a_k} -{\bf x}^{\bf u}\}$ be a generating set of $\mathcal{C}_{A}$ and let $n_{i}<n_{j}$. Then $C$ has Cohen-Macaulay tangent cone at the origin if and only if for every monomial $M=x_{i}^{v_{i}}x_{j}^{v_{^j}}x_{k}^{v_{k}}$, where $v_{i}<a_{i}$ and $v_{j}<a_{j}$, with $v_{i}n_{i}+v_{j}n_{j}+v_{k}n_{k} \in n_{1}+\mathcal{S}$ there exists a monomial $N$ with $1 \in {\rm supp}(N)$ such that $M-N \in I(C)$ and also ${\rm deg}(M) \leq {\rm deg}(N)$.

\end{theorem}

\begin{proposition} Suppose that $I(C)$ is given as in case 3. Let $\{x_1^{a_1} -
x_i^{a_i}, x_i^{a_i} - x_j^{a_j}, x_k^{a_k} -{\bf x}^{\bf u}\}$ be a generating set of $\mathcal{C}_{A}$ and let $n_{i}<n_{j}$. Assume that ${\rm deg}({\bf x}^{\bf u}) \leq a_{k}$ or/and $1$ belongs to the support of ${\bf x}^{\bf u}$. If \begin{enumerate} \item for every binomial $f=M-N \in \mathcal{I}$ with $1 \in {\rm supp}(M)$ we have that ${\rm deg}(N) \leq {\rm deg}(M)$ and \item for every monomial $x_{i}^{v_{i}}x_{j}^{v_{^j}}x_{k}^{v_{k}}$, where $v_{i}<a_{i}$, $v_{j}<a_{j}$ and $v_{k} \geq a_{k}$, with $v_{i}n_{i}+v_{j}n_{j}+v_{k}n_{k} \in n_{1}+\mathcal{S}$ there exists a monomial $N$ with $1 \in {\rm supp}(N)$ such that $M-N \in I(C)$ and also ${\rm deg}(M) \leq {\rm deg}(N)$,
\end{enumerate}
then $C$ has Cohen-Macaulay tangent cone at the origin.
\end{proposition}

\begin{corollary} Suppose that $I(C)$ is given as in case 3. Let $\{x_1^{a_1} -
x_i^{a_i}, x_i^{a_i} - x_j^{a_j}, x_k^{a_k} -{\bf x}^{\bf u}\}$ be a generating set of $\mathcal{C}_{A}$ and let $n_{i}<n_{j}$. Assume that $1$ belongs to the support of ${\bf x}^{\bf u}$. If \begin{enumerate} \item for every binomial $f=M-N \in \mathcal{I}$ with $1 \in {\rm supp}(M)$ we have that ${\rm deg}(N) \leq {\rm deg}(M)$ and \item $a_{k} \leq {\rm deg}({\bf x}^{\bf u})$,
\end{enumerate}
then $C$ has Cohen-Macaulay tangent cone at the origin.
\end{corollary}

\begin{theorem} Suppose that $I(C)$ is given as in case 4(a). Let $\{x_1^{a_1}-x_i^{a_i}, x_i^{a_i}-{\bf x}^{\bf u}, x_j^{a_j}-{\bf x}^{\bf v}, x_k^{a_k}-{\bf x}^{\bf w}\}$ be a generating set of $\mathcal{C}_{A}$. Then $C$ has Cohen-Macaulay tangent cone at the origin if and only if for every monomial $M=x_{i}^{z_{i}}x_{j}^{z_{^j}}x_{k}^{z_{k}}$, where $z_{i}<a_{i}$, with $z_{i}n_{i}+z_{j}n_{j}+z_{k}n_{k} \in n_{1}+\mathcal{S}$ there exists a monomial $N$ with $1 \in {\rm supp}(N)$ such that $M-N \in I(C)$ and also ${\rm deg}(M) \leq {\rm deg}(N)$.
\end{theorem}

\begin{proposition} Suppose that $I(C)$ is given as in case 4(a). Let $\{x_1^{a_1}-x_2^{a_2}, x_2^{a_2}-x_{3}^{u_{3}}x_{4}^{u_{4}}, x_3^{a_3}-{\bf x}^{\bf v}, x_4^{a_4}-{\bf x}^{\bf w}\}$ be a generating set of $\mathcal{C}_{A}$. Assume that $1 \in {\rm supp}({\bf x}^{\bf v})$ and $1 \in {\rm supp}({\bf x}^{\bf w})$. If \begin{enumerate} \item $a_{3} \leq {\rm deg}({\bf x}^{\bf v})$, \item $a_{4} \leq {\rm deg}({\bf x}^{\bf w})$ and \item for every binomial $f=M-N \in \mathcal{I}$ with $1 \in {\rm supp}(M)$ we have that ${\rm deg}(N) \leq {\rm deg}(M)$,
\end{enumerate}
then $C$ has Cohen-Macaulay tangent cone at the origin.
\end{proposition}

\begin{theorem} Suppose that $I(C)$ is given as in case 4(b). Let $\{x_1^{a_1} -
x_i^{a_i}, x_j^{a_j} - {\bf x}^{\bf u},  x_k^{a_k} - {\bf x}^{\bf v}\}$ be a generating set of $\mathcal{C}_{A}$. Then $C$ has Cohen-Macaulay tangent cone at the origin if and only if for every monomial $x_{i}^{z_{i}}x_{j}^{z_{^j}}x_{k}^{z_{k}}$, where $z_{i}<a_{i}$, with $z_{i}n_{i}+z_{j}n_{j}+z_{k}n_{k} \in n_{1}+\mathcal{S}$ there exists a monomial $N$ with $1 \in {\rm supp}(N)$ such that $M-N \in I(C)$ and also ${\rm deg}(M) \leq {\rm deg}(N)$.

\end{theorem}

\begin{proposition} Suppose that $I(C)$ is given as in case 4(b). Let $\{x_1^{a_1} -
x_i^{a_i}, x_j^{a_j} - {\bf x}^{\bf u},  x_k^{a_k} - {\bf x}^{\bf v}\}$ be a generating set of $\mathcal{C}_{A}$. Assume that $1 \in {\rm supp}({\bf x}^{\bf u})$ and $1 \in {\rm supp}({\bf x}^{\bf v})$. If \begin{enumerate} \item $a_{j} \leq {\rm deg}({\bf x}^{\bf u})$, \item $a_{k} \leq {\rm deg}({\bf x}^{\bf v})$, \item for every binomial $f=M-N \in \mathcal{I}$ with $1 \in {\rm supp}(M)$ we have that ${\rm deg}(N) \leq {\rm deg}(M)$ and \item for every binomial $f=M-N \in \mathcal{R}$ with $1 \in {\rm supp}(M)$ there is a monomial $P$ with $1 \in {\rm supp}(P)$ such that ${\rm deg}_{\mathcal{S}}(N)={\rm deg}_{\mathcal{S}}(P)$ and ${\rm deg}(N) \leq {\rm deg}(P)$,
\end{enumerate}
then $C$ has Cohen-Macaulay tangent cone at the origin.
\end{proposition}



\begin{thebibliography}{12}
	
	
	
	\bibitem{Ar}
	\textsc{F. Arslan,}
	\newblock \emph{Cohen-Macaulayness of tangent cones},
	\newblock Proc. Amer. Math. Soc. {\bf 128} (2000) 2243-2251
	
	\bibitem{ArMe}
	\textsc{F. Arslan, P. Mete,}
	\newblock \emph{Hilbert functions of Gorenstein monomial curves},
	\newblock Proc. Amer. Math. Soc. {\bf 135} (2007) 1993-2002.
	
	\bibitem{bayer}
	\textsc{D. Bayer, M. Stillman,}
	\newblock \emph{ Computation of Hilbert functions},
	\newblock J. Symbolic Comput. {\bf 14} (1992) 31-50.
	
	\bibitem{Bresinsky75}
	\textsc{H. Bresinsky,}
	\newblock \emph{Symmetric semigroups of integers generated by 4 elements},
	\newblock Manuscripta Math. \textbf{17}  (1975), no. 3, 205-219.
	
	\bibitem{Cas}
	\textsc{P. Pis\'{o}n Casares,}
	\newblock \emph{The short resolution of a lattice ideal},
	\newblock Proc. Amer. Math. Soc. \textbf{131} (2003) 1081-1091.
	
	\bibitem{Coc}
	\textsc{CoCoATeam,}
	\newblock \emph{CoCoA: A system for doing computations in commutative algebra},
	\newblock available at http://cocoa.dima.unige.it.
	
	
	\bibitem{Ga}
	\textsc{A. Garcia,}
	\newblock \emph{Cohen-Macaulayness of the associated graded of a semigroup rings},
	\newblock Comm. Algebra \textbf{10} (1982) 393-415.
	
	\bibitem{Gi}
	\textsc{P. Gimenez, H. Srinivasan,}
	\newblock \emph{A note on Gorenstein monomial curves},
	\newblock  Bull. Braz. Math. Soc. (N.S.) \textbf{4} (2014) 671-678.
	
	\bibitem{Herzog70}
	\textsc{J. Herzog,}
	\newblock \emph{Generators and relations of abelian semigroups and semigroup rings},
	\newblock Manuscripta Math. \textbf{3} (1970) 175-193.
	
	\bibitem{Her1}
	\textsc{J. Herzog,}
	\newblock \emph{When is a regular sequence super regular?}
	\newblock Nagoya Math. J. \textbf{83} (1981) 183-195.
	
	\bibitem{Her2}
	\textsc{J. Herzog, D. I. Stamate,}
	\newblock \emph{On the defining equations of the tangent cone of a numerical semigroup ring},
	\newblock J. Algebra \textbf{418} (2014) 8-28.
	
	\bibitem{Ja}
	\textsc{R. Jafari, S. Zarzuela,}
	\newblock \emph{On monomial curves obtained by gluing},
	\newblock Semigroup Forum \textbf{2} (2014) 397-416.
	
	\bibitem{KaOj}
	\textsc{A. Katsabekis, I. Ojeda,}
	\newblock \emph{An indispensable classification of monomial curves in $\mathbb{A}^{4}(K)$},
	\newblock Pac. J. Math. \textbf{268} (2014) 95-116.
	
	\bibitem{Ku}
	\textsc{E. Kunz,}
	\newblock \emph{The value semigroup of one dimensional Gorenstein ring},
	\newblock Proc. Amer. Math. Soc. \textbf{25} (1970) 748-751.
	
	\bibitem{Mo2}
	\textsc{S. Molinelli, D.P. Patil, G. Tamone,}
	\newblock \emph{On the Cohen-Macaulayness of the associated graded rings of monomial curves},
	\newblock Beitr\"age Algebra Geom.  \textbf{39} (1998) 433-446.
	
	\bibitem{Mo1}
	\textsc{S. Molinelli, G. Tamone,}
	\newblock \emph{On the Hilbert function of certain rings of monomial curves},
	\newblock J. Pure Appl. Algebra \textbf{101} (1995) 191-206.
	
	\bibitem{OST}
	\textsc{A. Oneto, F. Strazzanti, G. Tamone,}
	\newblock \emph{One-dimensional Gorenstein local rings with decreasing Hilbert function},
	\newblock arXiv: 1602.00334.
	
	
	\bibitem{Pa}
	\textsc{D.P. Patil, G. Tamone,}
	\newblock \emph{CM defect and Hilbert functions of monomial curves},
	\newblock Journal of Pure and Applied Algebra 215 \textbf{7} (2011) 1539-1551.
	
	\bibitem{Ro}
	\textsc{L. Robbiano, G. Valla,}
	\newblock \emph{On the equations defining tangent cones},
	\newblock Math. Proc. Cambridge Philos. Soc. \textbf{88} (1980) 281-297.
	
	\bibitem{Rossi}
	\textsc{M. Rossi,}
	\newblock \emph{Hilbert functions of Cohen-Macaulay local rings},
	\newblock Commutative Algebra and its Connections to Geometry, Contemporary Math \textbf{555} (2011), AMS, 173-200.
	
	\bibitem{Sa}
	\textsc{V. A. Sapko,}
	\newblock \emph{Associated graded rings of numerical semigroup rings},
	\newblock Comm. Algebra \textbf{10} (2001) 4759-4773.
	
	\bibitem{Sh}
	\textsc{L. Sharifan, R. Zaare-Nahandi,}
	\newblock \emph{Minimal free resolution of the associated graded ring of monomial curves of generalized arithmetic sequences},
	\newblock J. Pure Appl. Algebra \textbf{213} (2009) 360-369.
	
	\bibitem{She}
	\textsc{Y. H. Shen,}
	\newblock \emph{Tangent cone of numerical semigroup rings of embedding dimension three},
	\newblock Comm. Algebra \textbf{39} (2011) 1922-1940.
	
	\bibitem{Sturmfels95}
	\textsc{B. Sturmfels,}
	\newblock \emph{Gr\"obner bases and convex polytopes},
	\newblock volume~8 of \emph{University  Lecture Series.}
	\newblock American Mathematical Society, Providence, RI, 1996.
	
\end{thebibliography}
\end{document}